\documentclass[12pt]{article}
\usepackage[utf8]{inputenc}
\usepackage[T1]{fontenc}
\usepackage[english]{babel}
\usepackage{amsfonts}
\usepackage{amsmath}
\usepackage{amssymb}
\usepackage{amsthm}
\usepackage{mathtools}
\usepackage{enumitem}
\usepackage[]{authblk}
\usepackage{bbm}
\usepackage{xparse}
\usepackage{tikz}
\usepackage{centernot}
\usepackage{accents}
\usepackage{linegoal}
\usepackage[numeric,initials,nobysame,msc-links,abbrev]{amsrefs}
\renewcommand{\eprint}[1]{\href{https://arxiv.org/abs/#1}{arXiv:#1}}
\newcommand{\pageafter}[1]{#1~pp.}
\BibSpec{article}{%
+{} {\PrintAuthors} {author}
+{,} { \textit} {title}
+{.} { } {part}
+{:} { \textit} {subtitle}
+{,} { \PrintContributions} {contribution}
+{.} { \PrintPartials} {partial}
+{,} { } {journal}
+{} { \textbf} {volume}
+{} { \PrintDatePV} {date}
+{,} { \issuetext} {number}
+{,} { \pageafter} {pages}
+{,} { } {status}
+{,} { \PrintDOI} {doi}
+{,} { available at \eprint} {eprint}
+{} { \parenthesize} {language}
+{} { \PrintTranslation} {translation}
+{;} { \PrintReprint} {reprint}
+{.} { } {note}
+{.} {} {transition}
+{} {\SentenceSpace \PrintReviews} {review}
}
\BibSpec{collection.article}{%
+{} {\PrintAuthors} {author}
+{,} { \textit} {title}
+{.} { } {part}
+{:} { \textit} {subtitle}
+{,} { \PrintContributions} {contribution}
+{,} { \PrintConference} {conference}
+{} {\PrintBook} {book}
+{,} { } {booktitle}
+{,} { \PrintDateB} {date}
+{,} { \pageafter} {pages}
+{,} { } {status}
+{,} { \PrintDOI} {doi}
+{,} { available at \eprint} {eprint}
+{} { \parenthesize} {language}
+{} { \PrintTranslation} {translation}
+{;} { \PrintReprint} {reprint}
+{.} { } {note}
+{.} {} {transition}
+{} {\SentenceSpace \PrintReviews} {review}
}

\usepackage{verbatim}
\usepackage{subcaption}
\usepackage{hyperref}
\usepackage[capitalise]{cleveref}

\allowdisplaybreaks

\usetikzlibrary{arrows}
\usetikzlibrary[patterns]
\usetikzlibrary{decorations.pathreplacing}
\usetikzlibrary{hobby}

\definecolor{ffqqqq}{rgb}{1,0,0}
\definecolor{qqffqq}{rgb}{0,1,0}
\definecolor{ffffff}{rgb}{1,1,1}
\colorlet{ColorGray}{gray!30}

\newtheorem{mainthm}{Theorem}
\crefname{mainthm}{Theorem}{Theorems}
\newtheorem{thm}{Theorem}
\crefname{thm}{Theorem}{Theorems}
\newtheorem{cor}[thm]{Corollary}
\crefname{cor}{Corollary}{Corollaries}

\crefname{lem}{Lemma}{Lemmas}
\newtheorem{prop}[thm]{Proposition}
\crefname{prop}{Proposition}{Propositions}

\crefname{conj}{Conjecture}{Conjectures}

\crefname{ques}{Question}{Questions}
\theoremstyle{definition}

\crefname{defn}{Definition}{Definitions}
\newtheorem{rem}[thm]{Remark}
\crefname{rem}{Remark}{Remarks}

\crefname{ex}{Example}{Examples}

\crefname{obs}{Observation}{Observations}

\crefname{claim}{Claim}{Claims}

\crefname{ass}{Assumption}{Assumptions}

\numberwithin{thm}{section}

\newcommand{\cE}{\ensuremath{\mathcal E}}

\newcommand{\cG}{\ensuremath{\mathcal G}}

\newcommand{\bbE}{{\ensuremath{\mathbb E}} }

\newcommand{\bbN}{{\ensuremath{\mathbb N}} }

\newcommand{\bbP}{{\ensuremath{\mathbb P}} }
\newcommand{\bbQ}{{\ensuremath{\mathbb Q}} }
\newcommand{\bbR}{{\ensuremath{\mathbb R}} }

\newcommand{\bbT}{{\ensuremath{\mathbb T}} }

\newcommand{\bbZ}{{\ensuremath{\mathbb Z}} }

\newcommand{\ba}{\ensuremath{\mathbf{a} }}
\newcommand{\bb}{\ensuremath{\mathbf{b} }}
\newcommand{\bc}{\ensuremath{\mathbf{c} }}
\newcommand{\bd}{\ensuremath{\mathbf{d} }}
\newcommand{\be}{\ensuremath{\mathbf{e} }}

\newcommand{\bo}{\ensuremath{\mathbf{o} }}

\newcommand{\bu}{\ensuremath{\mathbf{u} }}
\newcommand{\bv}{\ensuremath{\mathbf{v} }}
\newcommand{\bw}{\ensuremath{\mathbf{w} }}
\newcommand{\bx}{\ensuremath{\mathbf{x} }}
\newcommand{\by}{\ensuremath{\mathbf{y} }}
\newcommand{\bz}{\ensuremath{\mathbf{z} }}
\newcommand{\bD}{\ensuremath{\mathbf{\D} }}

\renewcommand{\a}{{\ensuremath{\alpha}}}
\renewcommand{\b}{{\ensuremath{\beta}}}
\renewcommand{\d}{{\ensuremath{\delta}}}
\newcommand{\D}{{\ensuremath{\Delta}}}
\newcommand{\e}{{\ensuremath{\varepsilon}}}
\newcommand{\f}{{\ensuremath{\varphi}}}

\newcommand{\g}{{\ensuremath{\gamma}}}

\newcommand{\h}{{\ensuremath{\eta}}}
\renewcommand{\k}{{\ensuremath{\kappa}}}

\newcommand{\m}{{\ensuremath{\mu}}}
\newcommand{\n}{{\ensuremath{\nu}}}
\renewcommand{\o}{{\ensuremath{\omega}}}
\renewcommand{\O}{{\ensuremath{\Omega}}}

\newcommand{\s}{{\ensuremath{\sigma}}}
\renewcommand{\t}{{\ensuremath{\tau}}}
\newcommand{\x}{{\ensuremath{\xi}}}

\newcommand{\z}{{\ensuremath{\zeta}}}

\newcommand{\<}{\langle}
\renewcommand{\>}{\rangle}

\renewcommand{\leq}{\leqslant}
\renewcommand{\geq}{\geqslant}
\renewcommand{\le}{\leqslant}
\renewcommand{\ge}{\geqslant}
\DeclareDocumentCommand \to { o o } {%
  \IfNoValueTF {#1} {\IfNoValueTF{#2}{\rightarrow}{\xrightarrow[#2]}}%
{\IfNoValueTF{#2}{\xrightarrow{#1}}{\xrightarrow[#2]{#1}}}}
\DeclareDocumentCommand \nto { o } {%
  \IfNoValueTF {#1} {{\ensuremath{\centernot{\rightarrow}}} }{{\ensuremath{\centernot{\xrightarrow{#1}}}}}%
}

\newcommand{\pc}{{p_{\mathrm{c}}}}

\newcommand{\pcop}{{p_{\mathrm{c}}^{\mathrm{OP}}}}

\newcommand{\Uc}{\accentset{\circ}{U}}

\title{Generalised oriented site percolation\footnote{This work was supported by ERC Starting Grant 680275 ``MALIG.''}}

\author{Ivailo Hartarsky\thanks{\textsf{hartarsky@ceremade.dauphine.fr}} } 
\author{R\'eka Szab\'o\thanks{\textsf{szabo@ceremade.dauphine.fr}}}
\affil{CEREMADE, CNRS, Universit\'e Paris-Dauphine, PSL University\protect\\Place du Mar\'echal de Lattre de Tassigny, 75016 Paris, France}
\date{\vspace{-0.25cm}\today}
\begin{document}
\maketitle
\vspace{-0.75cm}

\begin{abstract}
We consider a generalised oriented site percolation model (GOSP) on $\bbZ^d$ with arbitrary neighbourhood. The key additional difficulties as compared to standard oriented percolation (OP) are the lack of symmetry and, in two dimensions, of planarity. We establish that, despite these deficiencies, in the supercritical regime GOSP behaves qualitatively like OP.
\end{abstract}

\section{Introduction}
\label{sec:intro}
Oriented percolation on $\bbZ^d$ is a classical model in probability theory and statistical physics, whose behaviour is relatively well understood with many of the main advances on the subject dating back to the 1980s (see \cites{Durrett84,Liggett05,Liggett99,Hinrichsen00} for comprehensive expositions). It is also essentially equivalent to the well-known contact process, but also linked to many other models and often used as a tool in proofs. 

In this work we study the supercritical phase of a natural generalisation of oriented site percolation on $\bbZ^d$ with arbitrary finite neighbourhood, which we define next. Our goal is to examine the importance of symmetry and planarity to the qualitative behaviour of oriented percolation. The generalisation is further motivated by its relations with probabilistic cellular automata and bootstrap percolation, as discussed in \cref{sec:background}.

\subsection{Model}
\label{subsec:model}
Our model of interest is \emph{generalised oriented site percolation} (GOSP) on $\bbZ^d$ for $d\ge 2$. The model is defined by a \emph{neighbourhood}---a finite set $X\subset \bbZ^d\setminus\{\bo\}$ ($\bo$ shall denote the origin of $\bbZ^d$) with $|X|\ge 2$ such that 
\begin{equation}
\label{eq:orientation}
    \exists \bu\in \bbR^d,\forall \bx\in X:\quad\<\bx,\bu\>>0,
\end{equation}
which ensures the orientation of the model, and a \emph{parameter} $p\in[0,1]$. For convenience we will always assume that $\bu=\be_d$, where $(\be_i)_{i=1}^d$ denotes the canonical basis of $\bbR^d$ and that the group generated by $X$ is $\bbZ^d$. This can be achieved by an invertible linear transformation of $\bbZ^d$ and, possibly, a restriction to a sublattice. We denote by $\bbP_p$ the product Bernoulli measure of parameter $p$ on $\bbZ^d$. The \emph{configuration} $\omega\in\Omega=\{0,1\}^{\bbZ^d}$ is assumed to be distributed according to this measure. We endow the vertex set $\bbZ^d$ with the locally finite translation-invariant oriented graph structure with edge set $\{(\ba,\ba+\bx):\ba\in\bbZ^d,\bx\in X\}$ generated by $X$ (see \cref{fig:example}). We refer to this graph as $\bbZ^d$ when $X$ is clear from the context and $\cG_{X}$ otherwise.

One can naturally identify $\o\in\O$ with the set of \emph{open sites} $\{\bx\in\bbZ^d:\omega_\bx=1\}\subset\bbZ^d$, all other sites being \emph{closed}. The open sites induce a subgraph of $\cG_X$ by keeping all edges between open sites. We can then introduce the following variant of the natural notion of being connected in this graph. For any $\ba, \bb\in \bbZ^d$ we say that \emph{$\ba$ infects $\bb$} (there is a path from $\ba$ to $\bb$) and write $\ba\to \bb$ for the event that there exists a sequence of open vertices $\ba_1,\dots, \ba_m=\bb$ such that $\ba_1-\ba\in X$ and $\ba_i-\ba_{i-1}\in X$ for all $i\in[2,m]$. Note that we do not require for $\ba$ to be open in order for $\ba\to \bb$ to occur. We make this choice so that $\ba\to \bb$ and $\bb\to \bc$ are independent for all $\ba,\bb,\bc\in\bbZ^d$.

For any $B\subset \bbZ^d$ we further define $\ba\to[B]\bb$ as $\ba\to \bb$ but with $\ba_i\in B$ for $i\in[1,m]$. We write $\ba\to[B]\infty$ for the existence of infinitely many $\bb$ such that $\ba\to[B]\bb$ and similarly for $\infty\to[B]\bb$. We further extend the notation by defining the event $C\to[B]D$ for $B,C,D\subset \bbZ^d$ as $\exists \bc\in C,\exists \bd\in D$ such that $\bc\to[B]\bd$. We say that $C$ \emph{percolates} in $B$ if $C\to[B]\infty$. We define the \emph{order parameter} \[\theta(p)=\bbP_p(\bo\to\infty),\]
the \emph{critical probability} \[\pc=\pc(X)=\inf\{p>0:\theta(p)>0\}\]
and say that \emph{there is percolation} at $p$ if $\theta(p)>0$ (by ergodicity this is equivalent to the a.s.\ existence of an infinite open path). Depending on the value of $p$, we may speak of \emph{subcritical}, \emph{critical} and \emph{supercritical} regimes. We focus on the study of the supercritical phase, where $\theta(p)>0$.

It is convenient to view the last coordinate of $\bbZ^d$ as the time in an interacting particle system. We therefore usually denote points in $\bbZ^d$ by $(\bx,t)$ with $\bx\in\bbZ^{d-1}$ and $t\in\bbZ$. Let $R=\max\{t\in\bbZ:(\bx,t)\in X\}$ be the \emph{range} of $X$. Consider the slab $S_t=\bbZ^{d-1}\times(\bbZ\cap [t,t+R))$ of width $R$ with normal vector $\be_d$, which we call \emph{time slab}, and denote $S=S_0$. Given an \emph{initial condition $A\subset S$} and a domain $B\subset\bbZ^d$, which we omit if $B=\bbZ^{d-1}\times[R,\infty)$, the \emph{state at time $t\in\bbN$} is 
\[_B\x^A_t=\left\{\bb\in S:\exists \ba\in A,\ba\to[B] \bb+t\be_d\right\},\]
so $(_{\bbZ^{d-1}\times[R,\infty)}\x^{A}_t)_{t=0}^\infty=(\xi^A_t)_{t=0}^\infty$ is a Markov chain with state space $\{0,1\}^S$. For simplicity if $A=\{\bo\}\subset S$, we write simply $\bo$ instead of $A$. Finally, in the supercritical phase it is useful to define
\[\bar\bbP_p=\bbP_p(\cdot|\forall t\ge 0,\x^\bo_t\neq\varnothing).\]

\subsection{Examples}
\label{subsec:examples}
Standard \emph{oriented percolation} in 2 dimensions (2dOP) can be defined by $X=\{(0,1),(1,1)\}$. However we will more customarily consider $X=\{(-1,1),(1,1)\}$ instead, which only spans half of $\bbZ^2$, but we will mostly disregard this minor detail. We denote by $\pcop$ the critical probability of 2dOP. In higher dimensions the situation is more ambiguous and at least the choices $X=\{\be_i:i\in\{1,\dots,d\}\}$; $X'=\{\be_d+\e \be_i:i\in\{1,\dots,d-1\},\e\in\{-1,1\}\}$ and $X''=X'\cup \{\be_d\}$ for the neighbourhood could be legitimately called \emph{$d$-dimensional oriented percolation} ($d$dOP). For concreteness, we will use $d$dOP to refer to $X''$ and simply OP for generic statements.

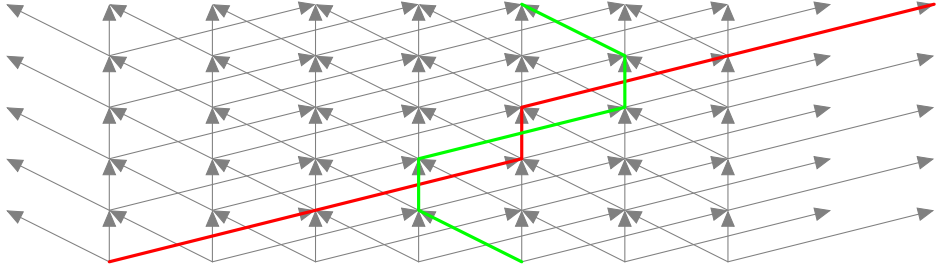
\begin{figure}
	\begin{center}
	\begin{tikzpicture}[line cap=round,line join=round,>=triangle 45,x=0.1\textwidth,y=0.05\textwidth]
\foreach \x in {1,...,7}
			{\foreach \y in {1,...,5}
			{\draw[->,color=gray] (\x,\y)--(\x-1,\y+1);
			\draw[->,color=gray] (\x,\y)--(\x,\y+1);
			\draw[->,color=gray] (\x,\y)--(\x+2,\y+1);}}
	\draw[very thick, color=red] (1,1) -- (5,3) -- (5,4) -- (7,5) -- (9,6);
	\draw[very thick, color=green] (5,1) -- (4,2) -- (4,3) -- (6,4) -- (6,5) -- (5,6);
	\end{tikzpicture}
	\caption{The graph $\cG_X$ on $\bbZ^2$ for $X=\{(-1,1),(0,1),(2,1)\}$. The two thickened paths cross although there is no common vertex nor an edge pointing from one to the other. Here $S=\bbZ\times\{0\}$, since $R=1$.\label{fig:example}}
	\end{center}
\end{figure}

As a prototype example of neighbourhood which is not covered by the classical approach, but handled here, we retain the two-dimensional GOSP defined by $X=\{(-1,1),(0,1),(2,1)\}$ (see \cref{fig:example}). It exhibits the two main additional difficulties of GOSP w.r.t.\ 2dOP: lack of symmetry w.r.t.\ the vertical axis and the non-planarity. The latter property is witnessed by the fact that paths may jump over each other without intersecting (see \cref{fig:example}).
\subsection{Results}
\label{subsec:results}
Denote by~$t^A(\bx)$ the \emph{hitting time} of~$\bx\in\mathbb{Z}^{d-1}$ from $A$:
\begin{equation}\label{eq:hittingtime}
  t^A(\bx):= \min\left\{t:(\bx, 0)\in\xi^A_t\right\},
\end{equation}
and define the following subsets of~$S$:
\begin{align}
H^A_t&{}:= \left\{(\bx, s)\in S: t^A(\bx)\leq t-s\right\},\label{eq:hitregion}\\
K^A_t&{}:= \left\{(\bx, s)\in S: \xi_t^A(\bx, s)=\xi_t^S(\bx, s)\right\},\label{eq:coupledregion}
\end{align}
which we refer to as \emph{hit} and \emph{coupled regions} with initial condition $A$ respectively. We omit $A$ if it is $\bo$. Our main result is the following.

\begin{mainthm}
\label{th:main}
Consider a GOSP in any dimension $d\ge 2$. For any $p>\pc$ there exists a deterministic convex compact set $U=U(p)\subset\mathbb R^{d-1}$ with non-empty interior such that for all~$\e>0$, $\bar\bbP_p$-a.s., for every~$t$ large enough it holds that
\begin{align}
\label{eq:main:lower}H_t\cap K_t&{}\supset(((1-\e)tU)\times[0,R))\cap\bbZ^d,\\
\label{eq:main:upper}\xi^\bo_t&{}\subset(((1+\e)tU)\times[0,R))\cap\bbZ^d.\end{align}
The function $p\mapsto U(p)$ is continuous on $(\pc,1]$ for the Hausdorff distance on non-empty compact subsets of $\bbR^{d-1}$. Furthermore, for any open set $O\subset U$, considering the cone  $C=\bigcup_{t>0}(tO\times\{t\})$, we have
\begin{equation}
\label{eq:restricted}\bbP_p\left(\exists \bx\in C,\bx\to[C]\infty\right)=1.
\end{equation}
\end{mainthm}

Our second result provides more precise information in the near-critical regime in two dimensions.
\begin{mainthm}
\label{th:2d}
For GOSP in two dimensions there exists $v\in\bbR$ such that 
\[\bigcap_{p>\pc} \Uc(p)=\{v\},\]
where $\Uc(p)$ is the interior of the limit shape from \cref{th:main}.
\end{mainthm}

\subsection{Organisation}
\label{subsec:organisation}
The paper is structured as follows. In 
\cref{sec:background} we provide the background for our work. In \cref{sec:preliminaries} we gather some preliminaries and notation. \Cref{sec:main,sec:2d} contain the proofs of \cref{th:main,th:2d} respectively.

The proofs are quite long and involve numerous intermediate results of independent interest. Inevitably, some of the steps are already known or require little or no new input as compared to existing arguments for OP or for the contact process. Nevertheless, we choose to also present these steps (without their proofs), so that the new ingredients we provide can be fitted into the global strategy and the reader is not obliged to scour the vast and entangled literature for all the ``well-known'' ingredients necessary. Moreover, in order not to disturb the flow of reasoning and to single out the novel contributions, we gather them in \cref{app}. Hence, specialists aware of classical results in two and more dimensions and of more recent developments around shape theorems may be able to directly consult \cref{app}.

\section{Background}
\label{sec:background}

We only discuss the supercritical regime, which is the focus of our work. Let us begin by emphasising that \cref{th:main,th:2d} are both known for OP and so are all intermediate results featuring in their proofs. More precisely, in the case of $d$dOP \cref{eq:main:lower,eq:main:upper} are due to \cites{Durrett83,Durrett91,Durrett82}; the continuity in \cref{th:main} was only recently established in \cite{Garet15}, based on \cites{Garet12,Garet14}; \cref{eq:restricted} was proved in \cite{Couronne04}*{Chapter 5}. Correspondingly, \cref{th:2d} for 2dOP was established in \cite{Durrett83} (see also \cites{Durrett84,Liggett05}).

Following progress on OP, natural generalisations similar to GOSP have often been considered. For the sake of comparability, in the present discussion, we focus on the most restrictive interesting case: GOSP with $X\subset\{(\ba,1):\ba\in \bbZ^{d-1}\}$, like the example of \cref{fig:example}. These models exhibit the main difficulties inherent to GOSP and are known as \emph{percolation probabilistic cellular automata} (PPCA), 2dOP being called \emph{Stavskaya's PCA} in this context \cites{Stavskaya71,Toom01, Toom95,Toom90,Taggi15,Taggi18}.

Bezuidenhout and Gray \cite{Bezuidenhout94} adapted the well-known renormalisation scheme of Bezuidenhout and Grimmett \cite{Bezuidenhout90} to show that in any dimension PPCA (and more general models) do not percolate at criticality. Their renormalisation will be the starting point of the proof of \cref{th:main}. In two dimensions an attempt at proving \cref{th:2d} and related results for PPCA (and more general models) was made by Durrett and Schonmann \cite{Durrett87}, themselves building on \cite{Durrett83}. Unfortunately, they imposed a restrictive technical assumption amounting to assuming that $X$ consists of consecutive sites. These neighbourhoods precisely lack the two main obstacles of GOSP---asymmetry and paths jumping over each other (see \cref{fig:example}). Furthermore, unaware of their work, Taggi \cites{Taggi15,Taggi18} claimed results for PPCA in two dimensions based on \cite{Durrett83}, as outlined in \cite{Durrett84}. Owing to non-planarity, his proof is only correct for neighbourhoods of consecutive sites. As we will see, \cite{Taggi15}*{Theorem 2.2} does indeed hold for all PPCA (and, more generally, GOSP), but requires a different treatment either based on higher dimensional techniques or on our enhancement of the approach of Durrett--Schonmann used to prove \cref{th:2d,th:main} respectively.

Let us note that GOSP are a particular case of \emph{boostrap percolation} \cites{Morris17,Hartarsky21,Schonmann92}. As established in \cite{Hartarsky21} (see particularly Remark 5.7 there), \cref{th:main,th:2d} on GOSP can be used to obtain results for more general bootstrap percolation models, particularly in conjunction with quantitative bounds on the limit shape $U$, as discussed in the first arXiv version of the present work \cite{Hartarsky21GOSParxiv}*{Section 5.7}. Other related models and generalisations of OP, to which much of the present approach applies can be found in \cites{Durrett82,Bezuidenhout94,Durrett80,Wierman85} (also see \cites{Deshayes15}).

\section{Preliminaries}
\label{sec:preliminaries}
\subsection{Duality}
\label{subsec:duality}
An important property of GOSP is that it is ``nearly'' self-dual (see \cites{Liggett05,Swart13} for background on duality). The dual of GOSP with neighbourhood $X$ can be thought of as a GOSP with paths moving ``backwards'' in time. More precisely, write $\ba\leadsto \bb$ if there exist $m\ge 0$ and $(\ba_i)_{i=0}^{m}$ with $\ba_0=\ba$ and $\ba_m=\bb$ such that for all $0\le i<m$ we have $\ba_i\in\o$ and $\ba_{i}-\ba_{i+1}\in X$. In other words, $\ba\to \bb$ iff $\bb\leadsto \ba$. Note that there are two differences with $\bb\to \ba$. Firstly, the steps are reversed: $\ba_{i+1}-\ba_{i}\in-X$. Secondly, for the dual connections we require that the initial site is open instead of the final one. Based on this notion we define the \emph{dual process} $(\tilde\x_t^A)$ again with state space $S$ but time coordinate $-\be_{d}$. We draw the reader's attention to the fact that this process does not have the same law as the primal process $(\x_t^A)$, for instance \[\tilde\theta(p)=\bbP_p(\bo\leadsto\infty)=p\bbP_p(\bo\to\infty)=p\theta(p).\]
However, up to such minor amendments all our results apply equally well to the dual process and we will use them as needed without systematically stating them.

\subsection{The contact process}
\label{subsec:discretisation}
OP is closely related to the \emph{contact process} (CP) \cites{Harris74,Liggett05,Liggett99}. The latter is often used to model epidemics on a graph: vertices are individuals, which can be healthy or infected. In this continuous time Markov dynamics infected individuals recover with rate~1 and transmit the infection to each neighbour with rate~$\lambda > 0$ (\emph{infection rate}). The CP admits a well-known graphical construction that is a space-time representation \cite{Liggett05}. We assign to each vertex~$\bv$ and ordered pair~$(\bu, \bv)$ of neighbours independent Poisson point processes~$D_\bv$ with rate~1 and~$D_{(\bu, \bv)}$ with rate~$\lambda$ respectively. For each atom~$t$ of~$D_\bv$ we place a ``recovery mark'' at~$(\bv, t)$ and for each atom of~$D_{(\bu, \bv)}$ we draw an ``infection arrow'' from~$(\bu, t)$ to~$(\bv, t)$. 
An \emph{infection path} is a connected path moving in the increasing time direction without crossing recovery marks, but possibly jumping along infection arrows in the direction of the arrow. Starting from a set of initially infected vertices~$A$, the set of \emph{infected} vertices at time~$t$ is the set of vertices~$\bv$ such that~$(\bv, t)$ can be reached by an infection path from some~$(\bu, 0)$ with~$\bu\in A$.

This representation can be thought of as a continuous time version of OP with infection paths in CP corresponding to paths in OP. Several of the results presented below are originally stated for CP but their proofs  transfer to discrete models with the following very minor adaptations.

Firstly, setting $\gamma=\max(\|\bx\|/t:(\bx,t)\in X)$, we clearly have that $\bo\to(\bx,t)$ implies $\|\bx\|\le \gamma t$, so, just like for CP, influence can spread at most linearly in time.

Secondly, since the group generated by $X$ is $\mathbb{Z}^d$, for all $n>0$ there exist a time $t$ and $\bv\in\bbZ^{d-1}$ such that
\begin{equation}\label{eq:l}
  \mathbb P_p\left(\xi_t^\bo\supset \bv+B_n\right)>0,
\end{equation}
where $B_n=([-n,n)^{d-1}\times[0,R))\cap\bbZ^d$.
This is the analogue of the fact that with positive probability the CP infects an arbitrarily large box in unit time.

Finally, for the CP one often needs to control the time an infection path spends at a vertex: either to ensure that it does not stay long at a vertex before jumping or that the path spends at least~$\d t$ time at a vertex during a time interval of length~$t$. The first assertion is trivial in discrete time as a path ``jumps immediately'' to the next vertex, and the discrete-analogue of the second assertion is visiting a vertex at least~$\lceil\d t\rceil$ times in a time interval of length~$t$.

\section{Proof of Theorem~\ref{th:main}}
\label{sec:main}

Throughout \cref{sec:main,sec:2d}, proofs will usually be omitted altogether when they only require minor changes (including the ones outlined in \cref{subsec:discretisation}). Nevertheless, we provide a sketch or at least a vague idea, whenever possible. The proofs requiring new ideas are gathered in \cref{app}.

The present section is structured as follows. In \cref{subsec:BG} we recall the Bezuidenhout--Grimmett renormalisation and its extension. We next derive several exponential bounds in \cref{subsec:restart} obtained based on restart arguments for later use. \Cref{subsec:shape} then completes the proof of the asymptotic shape result and its continuity from \cref{th:main}. Finally, \cref{subsec:cluster:properties} puts together all of the above with some further large deviation results to prove the percolation in restricted regions of \cref{th:main}. New ingredients needed in \cref{subsec:cluster:properties,subsec:shape,subsec:restart} are deferred to \cref{app:coupling:preliminary,app:tilting,app:LD}.

\subsection{Bezuidenhout--Grimmett renormalisation}
\label{subsec:BG}

We begin the study of the supercritical phase by briefly describing the well-known Bezuidenhout--Grimmett (BG) renormalisation. It was first introduced in~\cite{Bezuidenhout90} for the CP on~$\bbZ^d$ (see also~\cite{Liggett99}*{Sec. I.2}) and later generalised for translation-invariant finite-range attractive\footnote{A spin system is \emph{attractive} if adding extra sites in the initial condition only makes more sites infected by it at any later time (see e.g.\ \cite{Liggett99}*{Sec. I.1.}), as in the case of GOSP.} spin systems on~$\bbZ^d$ by Bezuidenhout and Gray \cite{Bezuidenhout94}, the latter reference being the most relevant for us. It is a construction that allows us to compare GOSP and 2dOP. The main idea is to show that GOSP when restricted to a sufficiently thick two-dimensional space-time slab dominates a supercritical 2dOP, which in turn implies that if there is percolation in 2dOP, then there is percolation in this restricted region. As all the results below are already known for 2dOP, this will entail numerous consequences for GOSP.

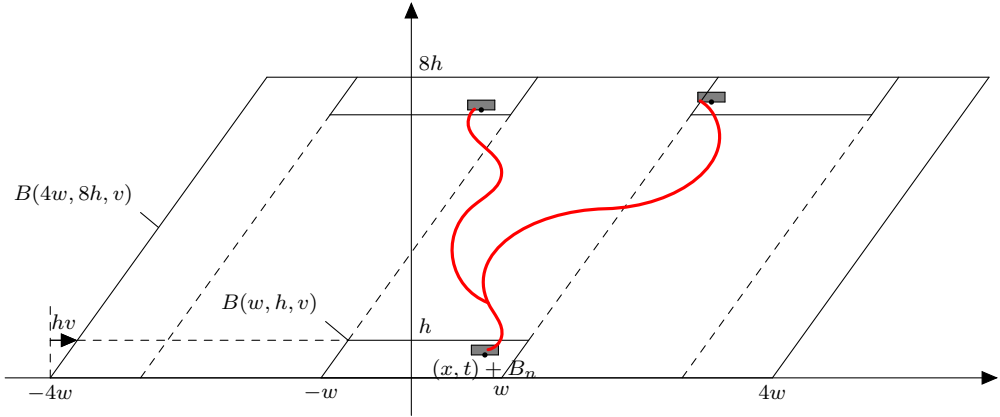
\begin{figure}
	\begin{center}
\begin{tikzpicture}[line cap=round,line join=round,>=triangle 45,x=0.6cm,y=0.5cm]
\draw[->] (-9,0) -- (13,0);
\draw[->] (0,-1) -- (0,10);
\fill[fill=black,fill opacity=0.5] (1.33,0.87) -- (1.33,0.61) -- (1.93,0.61) -- (1.93,0.87) -- cycle;
\fill[fill=black,fill opacity=0.5] (1.25,7.39) -- (1.25,7.13) -- (1.85,7.13) -- (1.85,7.39) -- cycle;
\fill[fill=black,fill opacity=0.5] (6.35,7.6) -- (6.35,7.34) -- (6.95,7.34) -- (6.95,7.6) -- cycle;
\draw (-2,0)-- (2,0);
\draw (-1.4,1)-- (2.6,1);
\draw (-8,0)-- (8,0);
\draw (-3.2,8)-- (12.8,8);
\draw (-3.2,8)-- (-8,0);
\draw (8,0)-- (12.8,8);
\draw (-1.4,1)-- (-2,0);
\draw (2,0)-- (2.6,1);
\draw (-1.2,8)-- (-1.8,7);
\draw (-1.8,7)-- (2.2,7);
\draw (2.8,8)-- (2.2,7);
\draw (6.2,7)-- (10.2,7);
\draw (10.2,7)-- (10.8,8);
\draw (6.2,7)-- (6.8,8);
\draw [->] (-8,1) -- (-7.4,1);
\draw (1.33,0.87)-- (1.33,0.61);
\draw (1.33,0.61)-- (1.93,0.61);
\draw (1.93,0.61)-- (1.93,0.87);
\draw (1.93,0.87)-- (1.33,0.87);
\draw (1.25,7.39)-- (1.25,7.13);
\draw (1.25,7.13)-- (1.85,7.13);
\draw (1.85,7.13)-- (1.85,7.39);
\draw (1.85,7.39)-- (1.25,7.39);
\draw (6.35,7.6)-- (6.35,7.34);
\draw (6.35,7.34)-- (6.95,7.34);
\draw (6.95,7.34)-- (6.95,7.6);
\draw (6.95,7.6)-- (6.35,7.6);
\draw [dashed] (6,0)-- (10.2,7);
\draw [dashed] (-6,0)-- (-1.8,7);
\draw [dashed] (-1.4,1)-- (2.2,7);
\draw [dashed] (6.2,7)-- (2.6,1);
\draw [dashed] (-8,0)-- (-8,2);
\draw [dashed] (-7.4,1)-- (-1.4,1);
\fill [color=black] (1.63,0.61) circle (1pt);
\fill [color=black] (1.55,7.13) circle (1pt);
\fill [color=black] (6.65,7.34) circle (1pt);
\draw [->] (-8,1) -- (-7.4,1);
\begin{scriptsize}
\draw (0,8) node [above right] {$8h$};
\draw (0,1) node [above right] {$h$};
\draw (8,0) node [below] {$4w$};
\draw (-8,0) node [below] {$-4w$};
\draw (2,0) node [below] {$w$};
\draw (-2,0) node [below] {$-w$};
\draw (-7.7,1.1) node [above] {$hv$};
\draw (-1.4,1)--(-1.9,1.5) node [above left] {$B(w,h,v)$};
\draw (-5.6,4)--(-6,4.4) node [above left] {$B(4w,8h,v)$};
\draw (1.63,0.7) node [below] {$(x,t)+B_n$};
\end{scriptsize}
\draw[color=red,very thick] (1.7,0.75) to [curve through={(2,1.1)..(1.7,2)..(4.3,4.5)..(6.7,7)}] (6.4,7.37);
\draw[color=red,very thick] (1.7,2) to [curve through={(1.3,4.5)..(2,5.5)..(1.3,7)}] (1.4,7.15);
\end{tikzpicture}
		\caption{The event described in \cref{thm:BGrenorm1} for~$d=2$. Note that in this case~$w$ and~$v$ are one-dimensional.} 
		\label{fig:BGrenorm}
	\end{center}
\end{figure}

Before proceeding to the renormalisation we will need a few geometric definitions. Our basic \emph{box} is $B_n=[-n, n)^{d-1}\times[0,R)$ for natural $n$, recalling $R$ from \cref{subsec:model} (although many of our regions will be defined as subsets of $\bbR^d$, we systematically refer to the integer points in them). For $\bw\in\mathbb Z^{d-1}$, $h\in\mathbb Z$ and~$\bv\in\bbR^{d-1}$ we further introduce the \emph{block} (see \cref{fig:BGrenorm})
\begin{equation}
B(\bw,h,\bv)=[0,h)(\bv,1)+\prod_{i=1}^{d-1}[-w_i,w_i)\times\{0\}\subset\bbR^d,
\label{eq:Bwhv}
\end{equation}
so that $B_n=B(\{n\}^{d-1},R,0)$. Note that if the model is symmetric we can always assume $\bv=0$. Here and below we write $w_i$ for the $i$\textsuperscript{th} coordinate of $\bw$ and similarly for other vectors.

The key theorem allowing comparison between the two models is as follows.
\begin{thm}\label{thm:BGrenorm1}
If $p$ is such that $\theta(p)>0$, then for all~$\e>0$ there exist positive integers $n,h$ and vectors $\bw\in\mathbb Z^{d-1}, \bv\in\mathbb R^{d-1}$ with~$n<w_{i}$ for all~$i$ and $h>R$ such that if~$(x, t)\in B(\bw, h, \bv)$, then
	\begin{multline}
	\bbP_p \big(\exists (\by, s)\in B(\bw, h, \bv)+7h(\bv, 1) \pm 2w_{d-1}\be_{d-1} \text{ such that }\\(\bx, t)+B_n \text{ infects }(\by, s)+B_n\text{ in }B(4\bw, 8h, \bv)\big)> 1-\e. \label{eq:BGrenorm1}
	\end{multline}
\end{thm}

In other words we can choose parameters such that  when considering the truncated process in~$B(4\bw, 8h, \bv)$ with high probability a box~$B_n$ centered at some~$(\bx, t)$ inside the block~$B(\bw, h, \bv)$ infects a copy of itself centered at either of the target blocks that are translates of~$B(\bw, h, \bv)$ (see \cref{fig:BGrenorm}). The proof of this theorem being quite long and technical, we direct the interested reader to \cite{Bezuidenhout94}, where it is established in a setting essentially including GOSP.

Recall that 2dOP is defined by $X=\{(-1,1),(1,1)\}$. We denote by $\z_k^{A}$ the set of (even) sites in $\bbZ^2$ with second coordinate $k$ infected by 2dOP with initial condition $A$. We are now ready for the 2dOP comparison of \cite{Bezuidenhout90}.

\begin{thm}[BG renormalisation]\label{thm:BGrenorm2}
Fix $q<1$ and assume that $p$ is such that \cref{eq:BGrenorm1} holds for $\e>0$ sufficiently small depending on $q$ and some $n,h,\bw,\bv$ as in \cref{thm:BGrenorm1}. Then the following holds for some $n,h,\bw,\bv$. For any initial condition $A\subset S$ we denote
\[A'=\left\{j\in2\bbZ:\exists (\bx,s)\in B(\bw,h,\bv)+2w_{d-1}j\be_{d-1}, (\bx,s)+B_n\subset A\right\}.\]
Then there exists a coupling of 2dOP $\z^{A'}$ of parameter $q$ and GOSP $\x^A$ such that for all $j\in\bbZ$ and $k$ \begin{multline*}
j\in \z^{A'}_k\text{ implies that }(\bx, 0)+B_n\subset\xi_t^A \\\text{ for some }(\bx,t)\in B(\bw,h,\bv)+7hk(\bv,1)+2w_{d-1}j\be_{d-1}.
\end{multline*}
In particular, $\theta(p)>0$.
\end{thm}

Informally, each site of 2dOP corresponds to a translate of the block $B(\bw, h, \bv)$ in GOSP. We can couple the two processes so that if a site in 2dOP is in the cluster of a vertex in~$\zeta^A$, then there is a box infected by~$A$ in the block corresponding to that site in the GOSP.

The proof of \cref{thm:BGrenorm2} is as in \cite{Liggett99}*{Theorem I.2.23} (also see \cite{Bezuidenhout90}). Indeed, one may construct the coupling by induction as follows. If $j\in\z_k$, then there is an infected copy of the box $B_n$ in the block corresponding to $j,k$, so we may apply the result of \cref{thm:BGrenorm1} to get that with probability $1-\e$ there will also be such infected boxes in the blocks corresponding to $j+1,k+1$ and $j-1,k+1$. It is easily seen (as the GOSP configuration is composed of independent variables) that the resulting process is a $1$-dependent 2dOP with parameter at least $1-\e$, so by a standard comparison between $1$-dependent and independent percolation \cite{Liggett97} we obtain \cref{thm:BGrenorm2} as desired.

It is useful to note that the BG renormalisation concerns only certain translates of the block $B(\bw,h,\bv)$. However, we may tile $\bbZ^d$ with disjoint blocks so that each tilted space-time slab of the form
\[\bigcup_{j,k\in\bbZ}B(\bw, h, \bv)+w_{d-1}j\be_{d-1}+kh(\bv, 1)\] is formed by 14 disjoint 2dOP lattices of blocks. We may perform the couplings of all the corresponding 2dOP processes with the same GOSP simultaneously as above so that sites in different 2dOP have a finite range dependence, hence they may be made independent by \cite{Liggett97}. In total, for $A\subset S$ we can couple $\x^A$ with independent 2dOP processes naturally indexed by $\bbZ^{d-2}\times \{1,\dots, 14\}$ with initial conditions corresponding to the parts of $A$ in each of the 2dOP lattices.

\subsection{Restart arguments}
Recall that $\x^A_t$ is the set of sites infected by $A$ at time $t$. The \emph{extinction time of $A$} is the absorption time of the chain started at $A$, that is \begin{equation}\label{eq:def:ta}
\t^A=\min\left\{t\geq 0:\x_t^A=\varnothing\right\}\in\{0,1,\dots\}\cup\{\infty\}.\end{equation}
We say that the process started from $A$ \emph{dies out} if $\t^A<\infty$ and \emph{survives} otherwise. 

\label{subsec:restart}The BG renormalisation allows us to use a so-called \emph{restart argument} that can be described as follows. We let GOSP $(\xi_t^A)$ evolve until we find an infected translate of the box~$B_n$ (which by~\cref{eq:l} has a strictly positive probability, thus it happens after at most a geometrically distributed number of steps) or the process dies out. If we infect a box, we start the 2dOP process~$(\zeta_k)$ of~\cref{thm:BGrenorm2} (with appropriate initial condition) from the corresponding block coupled with~$(\xi^A_t)$. If~$(\zeta_k)$ percolates, then~$(\xi^A_t)$ percolates as well. If $(\z_k)$ dies out and~$(\xi^A_t)$ still survives, we restart the procedure. We repeat this until either the GOSP dies out or the renormalisation yields a percolating 2dOP (since the parameter $q$ of $(\zeta_k)$ is supercritical, $q>\pcop$, this will happen after at most a geometric number of trials).

We will use this technique to transfer properties from 2dOP to GOSP. The exponential bounds we prove next in this section (like all other results) are already established for 2dOP \cite{Durrett83} and the $d$-dimensional CP \cite{Durrett91}. Recall that~$\tau^A$ is the extinction time of the set~$A$.

\begin{thm}[Exponential death bounds]
\label{thm:GOSPproperties}
For every~$p>\pc$ there exists a constant~$\e=\e(p)>0$ such that for all $A\subset S$ and $t\ge R$ it holds that
	\begin{align}\label{eq:property1}
\mathbb{P}_p\left(t\le \tau^A<\infty\right)&{}\leq e^{-\e t},
\\\label{eq:property2}
	    \mathbb{P}_p\left(\tau^A<\infty\right)&{}\leq e^{-\e |A|}.
	\end{align}
\end{thm}

The proof of \cref{thm:GOSPproperties} goes along the same lines as the proof of \cite{Liggett99}*{Theorem I.2.30}, using the restart argument. For~\cref{eq:property1}, on~$\{\tau^A<\infty\}$, we can bound~$\tau^A$ by the sum of the number of steps until we find an infected box and the survival time of the coupled 2dOP in each round. As 2dOP satisfies~\cref{eq:property1} and the required quantities of the restart argument are bounded by geometric random variables, we get the desired exponential decay.

For~\cref{eq:property2} first tile~$\bbZ^d$ with the disjoint translates of the space-time slab
\[T=\bigcup_{j,k\in\bbZ}B(4\bw, 7h, \bv)+4w_{d-1}j\be_{d-1}+7hk(\bv, 1)\]
and consider the processes restricted to these slabs with initial conditions corresponding to the parts of~$A$ in each slab. Observe that these processes are independent and~$\tau^A$ can be bounded from below by the maximum of their extinction times. Therefore, it is enough to show the analogue of~\cref{eq:property2} for~$A\subset T$.
As in~\cref{eq:l} we can show that there exists a~$t$ depending on $\bw$ and $h$, but not on $A$, such that with strictly positive probability every vertex in~$A\subset T$ can infect a box~$B_n$ in
\[t(\bv,1)+\bigcup_{j\in\bbZ}B(\bw, h, \bv)+w_{d-1}j\be_{d-1}.\]
Thus, (save for an exponentially unlikely event) for some $\d>0$ at least~$\delta|A|$ disjoint blocks at the same time contain an infected box. We start the 2dOP process
of \cref{thm:BGrenorm2} from all these blocks. Observe that~$\{\tau^A<\infty\}$ can only happen if all the coupled 2dOP process dies, but since~\cref{eq:property2} holds for 2dOP, this has exponentially small probability in the size of~$A$.
 
\begin{rem}
\Cref{eq:property1} implies that the law of $\x_t^S$ converges to the upper invariant measure of the process (corresponding to the distribution of sites $\bx\in S$ such that $\bx\leadsto\infty$) exponentially fast in $t$.
\end{rem}

The next result is the analogue of condition (a) of Lemma 5.1 in \cite{Durrett91}.
\begin{thm}\label{thm:GOSPproperties2}
Let~$\x$ and $\tilde\x$ be independent primal and dual GOSP. Then for every~$p>\pc$ there exist constants~$\e, c, C>0$ and a vector $\bv\in\bbR^{d-1}$ depending on $p$ such that for all integer~$t>0$ and $A, B\subset S$ satisfying $\max_{\ba\in A, \bb\in B}\| \ba-\bb\|<ct$ we have
\begin{equation}
\label{eq:property4bis}
\bbP_p\left(\xi_t^A\cap\tilde{\xi}_t^{B+(2t\bv+\bz_t,0)}=\varnothing, \xi_t^A\neq\varnothing, \tilde{\xi}_t^{B+(2t\bv+\bz_t,0)}\neq\varnothing\right)\leq Ce^{-\e t},
\end{equation}
where $\bz_t\in\bbR^{d-1}\times\{0\}$ is such that $2t\bv+\bz_t\in\bbZ^{d-1}$ and $\|\bz_t\|^2\le (d-1)/4$.
\end{thm}
It is important to note that due to the lack of symmetry this result is more technical for GOSP than for $d$dOP. We leave the proof to \cref{app:coupling:preliminary} and only indicate that it relies mainly on the BG renormalisation, restart argument, \cref{eq:property1} and several properties known for 2dOP, which will be established below for GOSP.

\begin{rem}
We can use \cref{thm:GOSPproperties2} to prove that the infinite cluster is unique. Together with $\theta(\pc)=0$ following from \cref{thm:BGrenorm1,thm:BGrenorm2}, this customarily yields that $\theta:p\mapsto\theta(p)$ is continuous on $[0,1]$. This was first established for $d$dOP in~\cite{Grimmett02}.
\end{rem}

Recall the hit and coupled regions of~\cref{eq:hitregion,eq:coupledregion}.
\begin{thm}[At least linear growth]\label{thm:GOSPproperties3}
For every~$p>\pc$ there exist constants~$\e, c>0$ and a vector $\bv\in\bbR^{d-1}$ depending on $p$ such that for all $t>0$ and $\bx\in\bbZ^{d-1}$ such that $\| \bx-t\bv  \| <ct$ it holds that
	\begin{align}
	    \mathbb{P}_p\left((\bx,0)\not\in H_t, \tau^\bo=\infty\right)&{}\leq e^{-\e t}, \label{eq:property3}\\
        \mathbb{P}_p\left((\bx,0)\not\in K_t, \tau^\bo=\infty\right)&{}\leq e^{-\e t}.  \label{eq:property5}
	\end{align}
\end{thm}
This result is also a consequence of \cref{thm:GOSPproperties2}. The proof is an adaptation of the proof of conditions~(c) and~(d) of Theorem 5.2 in~\cite{Durrett91} for the CP (see also \cite{Durrett82}). For~\cref{eq:property3} recall the definition of the dual process~$\tilde\xi$ from \cref{subsec:duality} and note that for any $0\le s\le t$ the event $\{\exists \by\in S,\bo\to \by+s\be_d,(\bx,t)\leadsto \by+s\be_d\}$ is equivalent to~$\bo\to(\bx,t)$. By a restart argument we can find a time $r$ such that $(\bx,r)\leadsto\infty$ and up to an exponentially unlikely event we can take $0\leq t-r<ct$. Then the conclusion follows from \cref{eq:property4bis} with $A=\{\bo\}$, $B=\{(\bx-2r\bv,0)\}$ and $t$ in the theorem equal to $r/2$.
 
For~\cref{eq:property5} observe that for $\o\subset\bbZ^{d-1}\times\{R,R+1,\dots\}$
\[\left\{(\bx,0)\not\in K_{2t}, \tau^\bo=\infty\right\}=\left\{\bo\nto(\bx, 2t), (\bx, 2t)\leadsto S, \bo\to\infty\right\}.\]
The conclusion then follows again from \cref{eq:property4bis} with $A=\{\bo\}$ and $B=\{(\bx-4t\bv,0)\}$.

Finally, for completeness let us mention a consequence of the restart arguments concerning GOSP on tori. Recall that we assume that the direction $\bu=\be_d$, and let $\bbT^{d-1}_n=(\bbZ/n\bbZ)^{d-1}$ denote the $(d-1)$-dimensional discrete torus of side $n$. Consider GOSP on the graph with vertex set $\bbT^{d-1}_n\times \bbZ$ obtained as the quotient of $\cG_X$. 
The \emph{extinction time} is defined as
\[\tau^{\bbT}=\sup\left\{t\ge 0: \bbT^{d-1}_n\times\{0,-1,-2,\dots\}\to\bbT^{d-1}_n\times\{t\}\right\}.\]

\begin{cor}
\label{th:tori}
For all $p<\pc$ there exists $c(p)$ such that
\begin{equation}
\label{eq:meta:subcrit}
\frac{\tau^\bbT}{\log n}\to[\bbP_p][n\to\infty]\frac{d-1}{c(p)},
\end{equation}
and for all~$p>\pc$ there exist~$c,C\in(0,\infty)$ such that
\begin{gather}
\label{eq:meta1}
\frac{\t^\bbT}{\bbE_p[\t^\bbT]}\to[(\mathrm{d})][n\to\infty]\cE,\\
e^{cn^{d-1}}<\bbE_p\left[\t^\bbT\right]<e^{Cn^{d-1}},\label{eq:meta2}
\end{gather}
for $n$ large enough, where $\cE$ is the standard exponential distribution.
\end{cor}

\Cref{eq:meta:subcrit} follows as in~\cite{Durrett88a} (see also \cite{Liggett99}*{Theorem I.3.3}) from the subcritical result established in~\cites{Aizenman87,Menshikov86}: for all~$p<p_c$ there exists~$c(p)$ such that
\begin{equation}
\label{eq:subcrit}
-\lim_{t\to\infty}\frac{1}{t}\log \bbP_p(\tau^\bo\ge t)=c(p)>0.
\end{equation}
\Cref{eq:meta1} was proved for the CP in two dimensions in \cite{Schonmann85} (see also \cite{Durrett88b}), while in $d$ dimensions this was done in \cite{Mountford93} (also see \cite{Simonis96} for subsequent development). \Cref{eq:meta2} was proved in \cite{Durrett88a} for the two-dimensional CP and in \cite{Chen94} in $d$ dimensions. The proofs rely on \cref{thm:GOSPproperties3,thm:GOSPproperties,thm:BGrenorm2} (see \cite{Hartarsky21GOSParxiv}*{App. A.4} for a sketch following \cite{Mountford93}).

Let us note that for the CP on a finite box (in our setting this corresponds to cutting the bonds crossing the boundary of a fundamental domain of the torus) in \cites{Durrett88b,Mountford99} it was established that in fact $\log\bbE_p[\t]/n^{d-1}$ converges as $n\to\infty$. However, for GOSP considering a box is either inappropriate or requires tilting the lattice first, making the result somewhat unnatural and unhandy due to the implicit definition of the tilting direction, which may even depend on $p$ as $p\to\pc$. It would appear that proving the existence of the above limit on the torus is unknown even for OP and CP.

\subsection{Asymptotic shape}
\label{subsec:shape}
With the results of~\cref{subsec:restart} at hand we are ready to prove the asymptotic shape theorem and the continuity of the limit shape, that is~\cref{eq:main:lower,eq:main:upper} and the continuity result of~\cref{th:main}. These results are known for the CP. However, certain issues arise due to the possibility that the model may have a ``drift,'' e.g.\ if the convex envelope of the neighbourhood $X$ does not intersect the line $\bbR \be_d$. This problem is absent if we can take $\bv=0$ in \cref{thm:GOSPproperties3}. For simplicity, in \cref{subsec:shape} we only briefly recall the arguments used to prove the desired results under this additional assumption, leaving out the minor changes described in \cref{subsec:discretisation}. Thus, we leave the new input needed for removing the assumption $\bv=0$ to \cref{app:tilting}.

It was proved in~\cite{Durrett82} for permanent one-site growth procesess (translation invariant, attractive processes with local rules, with~$\varnothing$ absorbing state and positive probability of survival) that the exponential estimates from~\cref{eq:property1,eq:property3,eq:property5}  with $\bv=0$ imply the shape theorem: \cref{eq:main:lower,eq:main:upper}. The idea is that, given these estimates, the hitting times are subadditive, stationary and integrable. Then, using subadditive ergodic theory~\cite{Kingman73}, one can prove that for $\bx\in\bbZ^{d-1}$
\begin{equation}\label{eq:mu}
\frac{t(n\bx)}n\to\mu(\bx) \quad\bar\bbP_p\text{-a.s.} \end{equation}
The \emph{time constant} $\m(\bx)$ can be extended into a norm on~$\mathbb{R}^{d-1}$ with unit ball $U$, yielding the result for the hit region. Then we can argue that there are a lot of vertices around the boundary of the cone defined by~$U$ that are reached from the origin and by~\cref{eq:property1} survive. Using~\cref{eq:property5} we can conclude that the union of the coupled regions of these vertices eventually covers~$(1-\e)tU$.

Our next goal is to prove that the limit shape $U$ is continuous in $p$. For this, we will require a quantity called essential hitting time. It was first introduced by Garet and Marchand in~\cite{Garet12}, inspired by \cite{Kuczek89}, to prove shape theorems in a more difficult setting. Using this notion, they later proved large deviation inequalities~\cite{Garet14} and continuity of the asymptotic shape~\cite{Garet15}. We next discuss these results still under the assumption that $\bv=0$ in \cref{thm:GOSPproperties3}.

Roughly speaking (see \cite{Garet12} for the correct definition), under $\bar\bbP_p$ the \emph{essential hitting time} $\s(\bx)$ of $\bx\in\bbZ^{d-1}$ is a time such that $\bo\to(\bx,\s(\bx))\to\infty$. Crucially, the essential hitting time is nearly subadditive \cite{Garet12}*{Theorem 2}.
Using this property, one can show that~$\bar\bbP_p$-a.s., as $n\to\infty$, $\sigma(n\bx)/n$ converges. Controlling the discrepancy between the essential hitting time and the hitting time \cite{Garet12}*{Proposition~17}, we can conclude that the limit is also the one of $t(n\bx)/n$. This control of~$\sigma(\bx)-t(\bx)$ further allows us to bound the moments of~$\sigma(\bx)$ under $\bar\bbP_p$ and to get exponential estimates for the essentially hit region analogous to~\cref{eq:property3} with $\bv=0$ \cite{Garet12}*{Corollary~20 and~21}.

Relying on the almost subadditivity of~$\sigma(\bx)$, one may establish large deviation results corresponding to \cite{Garet14}*{Theorems 1.1 and 1.4}, still under the assumption $\bv=0$ to be removed in \cref{app:tilting}.
\begin{thm}
\label{thm:largedev}
For every~$p>\pc$ and every~$\e>0$ there exist constants $c, C>0$ such that for any~$\bx\in\mathbb{Z}^{d-1}$ and $t\ge 1$
\begin{align*}
\bar\bbP_p\left(K_t\cap H_t\supset (((1-\e)tU)\times [0,R))\cap\bbZ^d\right)&{}\geq 1-Ce^{-ct},\\
\bbP_p\left(\x^\bo_t\subset (((1+\e)tU)\times [0,R))\cap\bbZ^d\right)&{}\geq 1-Ce^{-ct},
\end{align*}
where $U$ is as in \cref{th:main}.
\end{thm}
Finally, one can show the continuity of the limit shape in~\cref{th:main} as in \cite{Garet15}*{Theorem 1}, recycling much of the proof of \cref{thm:largedev}.

\subsection{Percolation in restricted regions}
\label{subsec:cluster:properties}

Relying on the results of~\cref{subsec:restart}, we next establish large deviations for the infinite cluster density, which we then use to prove~\cref{eq:restricted} of~\cref{th:main}.

\begin{thm}
Let
\begin{equation}
\label{eq:def:Yn}
Y_n:=\left|\left\{(\bx,t)\in B_n:(\bx,t)\leadsto\infty\right\}\right|/|B_n|.
\end{equation}
\label{thm:LD}
For all $p>\pc$ there exists a convex function $\f:[0,1]\to{} [0,\infty)$ such that $\f(a)=0$ if and only if $a=p\theta(p)=\tilde\theta(p)$ and for all $a<b$ in $[0,1]$,
\[\lim_{n\to\infty}\frac{\log\bbP_p(Y_n\in[a, b])}{n^{d-1}}=-\inf_{x\in[a,b]}\varphi(x).\]
\end{thm}
Most of this result is very general and holds for any translation invariant attractive spin system, as established in \cite{Lebowitz88}. Roughly speaking, the existence of the limit follows from the fact that if several boxes have $Y_n\ge a$, then so does their union; the convexity follows similarly, asking for one box with $Y_n\ge x$ and one with $Y_n\ge y$ and considering their union. The relevance of $\tilde\theta(p)$ comes from cutting a large box into smaller ones and using attractiveness to replace boundary conditions by maximal ones, in order to enable the use of large deviations results for i.i.d.\ random variables. Indeed, the invariant measure of GOSP in a large box with infected boundary condition still infects sites ``far from the boundary'' with probability close to $\tilde\theta(p)$. This follows from the fact that the upper invariant measure of the infinite volume attractive process must dominate the (decreasing) limit of these invariant measures as the size of the box diverges (see \cite{Liggett05}*{Theorem III.2.7}).

The only somewhat model-specific property is the fact that for any $a<\tilde\theta(p)$ there exists $c>0$ such that for all $n$ we have
\begin{equation}
\label{eq:LD:reduced}
\bbP_p(Y_n\le a)\le e^{-cn^{d-1}}.
\end{equation}
This was established for 2dOP in \cite{Durrett88}. Unfortunately, the argument is 2-dimensional, so we provide a proof for GOSP in any dimension (and in particular $d$dOP), which appears to be novel. This is done via a new renormalisation in \cref{app:LD}, relying on \cref{thm:GOSPproperties,thm:GOSPproperties3}.

\begin{rem}One can further study fluctuations of the density. For translation invariant attractive spin systems on~$\bbZ^d$ \cite{Lebowitz87} examined when, starting from the upper invariant measure, we reach a value of $Y_n$ smaller than $\tilde\theta(p)$. Later this result was extended to upper fluctuations in~\cite{Galves89} for the CP. These proofs rely on \cref{thm:LD} and can be adapted to GOSP. We also direct the reader to \cite{Couronne04}*{Chapter 5} for information regarding the properties and shape of large finite clusters and more large deviations.\end{rem}

Now we are ready to prove a more geometric property of the infinite cluster,~\cref{eq:restricted} of~\cref{th:main}, establishing that percolation occurs in restricted regions. This result was proved for $d$dOP in \cite{Couronne04}*{Theorem 1.3 of Chapter 5}, but given the results available to us, we may directly retrieve it (for GOSP) from \cref{thm:GOSPproperties,thm:largedev,thm:LD} as follows. Fixing some $\bu\in O$ and $\d>0$ small, for a site $(\bx,t)$ at distance at most $\d^2 t$ from $(t\bu,t)$ surviving for time $\d t$, by \cref{eq:property1,thm:largedev} it is likely that its coupled region contains a box of side $\d^2 t$ centered at $((1+\d)t\bu,(1+\d)t)$. By \cref{thm:LD}, it is likely that at least $\d^{2d}t^{d-1}$ sites in that box are infected by $(\bx,t)$ and, since $\d$ is small and $(\bx,t)$ is at distance of order $t$ from the boundary of $C$, this has to happen inside $C$. Finally, by \cref{eq:property2}, some of those sites is likely to survive. Repeating this procedure to infinity and recalling that the probability of failing at each step is exponentially small in $t$, we obtain \cref{eq:restricted} of \cref{th:main}, as desired.

\section{Proof of Theorem~\ref{th:2d}}
\label{sec:2d}
In this section we assume $d=2$. In this case one can say more about GOSP based on techniques for 2dOP, for which \cite{Liggett05}*{Chapter VI} and \cite{Durrett84} are excellent references. In~\cref{sec:edge} we gather some standard preliminaries. In~\cref{subsec:DS} we introduce an alternative renormalisation technique, whose refinement enables us to prove~\cref{th:2d}.

\subsection{Edge speed}
\label{sec:edge}
Define the \emph{right edge} of the process as
\[r_t=\max\left\{x\in \bbZ: \exists y\in \{0,\dots,R-1\}, (x,y)\in\x_t^{S^-}\right\},\]
where $S^-=((-\infty,0]\times[0,R))\cap\bbZ^2$ is the left half of $S$. Similarly define the \emph{left edge} $l_t$ as the minimum of $x$ above with $S^+=([0,\infty)\times[0,R))\cap\bbZ^2$ instead of $S^-$. It is important to note that, as discussed in \cref{subsec:examples}, although the model is two-dimensional, it is \emph{not} planar and paths may jump over each other without crossing (recall \cref{fig:example}). Nevertheless, the right and left edges do retain some of their properties from the 2dOP case.

\begin{thm}[Edge speed]
\label{th:alpha}
For any $p\in[0,1]$ there exists
\[\a=\lim_{t\to\infty}\frac{\bbE_p[r_t]}{t}=\inf_{t\ge 1}\frac{\bbE_p[r_t]}{t}\in[-\infty,\infty).\]
Moreover, $r_t/t\to[t\to\infty]\a$ a.s.\ and if $\a>-\infty$, then $\bbE_p\left[\left|\frac{r_t}{t}-\a\right|\right]\to[t\to\infty]0$. Similar statements hold for $\beta=\lim\bbE_p[l_t]/t$.
\end{thm}
The proof (and statement) is identical to \cite{Liggett05}*{Theorem VI.2.19} and is a consequence of a subadditive ergodic theorem due to Durrett \cite{Durrett80} (see particularly Theorem 6.1 thereof). The idea is to introduce a version of the right edge between time $s$ and $t$ which, contrary to $r_t$, is subadditive in an appropriate sense.

We next show that the two-dimensional approach coincides with the more general one from the previous section.
\begin{thm}
\label{th:identification}
For any $p>\pc$ the limit shape $U$ from \cref{th:main} and the edge speeds $\a,\b$ from \cref{th:alpha} satisfy $U=[\b,\a]$.
\end{thm}
To see this, note that by \cref{thm:GOSPproperties3} with positive probability $\bo\to\infty$ and at all times the coupled region is large enough to ensure that the right and left edges are infected by $\bo$. We can then conclude, since \cref{th:main,th:alpha} are almost sure statements.

A notable advantage of having the edge representation of the limit shape is the following result.
\begin{thm}
\label{th:strict:monotonicity}The right edge speed $\a$ is strictly increasing on $(\pc,1)$.
\end{thm}
The proof is very similar to \cite{Durrett84}*{Eq.~(12)} and was reiterated in \cite{Durrett87} in a setting including PPCA. It proceeds in two steps. Firstly, one shows by a clever but simple algebraic manipulation that adding a vertical column of sites to any initial condition entirely on its right increases $\bbE_p[r_t]$ by at least $1$ for all $t$. This property only relies on the fact that the process is additive in the sense that $\x^{A\cup B}_t=\x_t^{A}\cup\x_t^B$ for all $A,B,t$. Secondly, one observes that if $p$ is increased by a small amount $\delta$, it may happen that the additional vertices opened by increasing it lead precisely to adding such a vertical column in $\x_t$ to the right of $r_t$ (corresponding to parameter $p$).

\subsection{Alternative renormalisation}
\label{subsec:DS}

In two dimensions it is possible to study the supercritical phase via a more elementary renormalisation scheme than the BG one. For 2dOP this approach due to Durrett and Griffeath \cite{Durrett83} is used classically to derive most of the results stated above in that setting. However, applying this renormalisation to GOSP (in two dimensions) turns out to be quite tricky. Let us first give the rough lines of the renormalisation before explaining what goes wrong for GOSP and how to address it.

For 2dOP, let us assume that $p$ satisfies $\a(p)>\b(p)$. By \cref{th:alpha} we have that $r_t/t\to\a$ a.s. Moreover, one can show (see \cites{Durrett83,Durrett84}) that for all $\e>0$ there exists $c>0$ such that for all $t>0$
\begin{equation}
\label{eq:edge:LD}
\bbP_p(r_t>(\a+\e)t)\le e^{-ct}.
\end{equation}
We may then establish (see \cref{fig:DS1}) that for $L$ large enough depending on $\e>0$, the box $B(\e L,L,\a)$ is crossed from bottom to top by an open path with high probability, namely for $\e>0$
\begin{equation}
\label{eq:crossing}
\lim_{L\to\infty}\bbP_p\left(_{B(\e L,L,\a)}\x^{S^-}_L=\varnothing\right)=0.
\end{equation}
Indeed, by \cref{eq:edge:LD}, it is forbidden for the right edge to leave the box on one side; by \cref{th:alpha} the right edge at time $L$ is likely to be in the middle of the top side of the box; while if the path reaching the right edge at time $L$ leaves the box on the other side, that would imply that the path necessarily went faster than allowed by \cref{eq:edge:LD}, in order to make up for the delay (the last idea is due to Gray \cite{Durrett87}).

Hence, we have that with probability close to $1$ long thin boxes with tilting $\a$ (and similarly for $\b$) are crossed. This reasoning is perfectly valid for GOSP. In order use such boxes to construct a renormalisation, one places around each renormalised vertex two of them directed by $\a$ and $\b$ and says the vertex is open if they are crossed by paths (see \cite{Durrett84}*{Fig. 7} or \cite{Durrett83}*{Fig. 1}). For 2dOP it is then clear that if the resulting renormalised 2dOP percolates, then so does the original one. Indeed, one can switch from the path in one box to another as soon as they intersect, which is necessarily the case for planar graphs such as the one associated to 2dOP.

It is not hard to see that the argument remains valid for PPCA with neighbourhood consisting of consecutive sites of the form $(x,1)$. However, for GOSP with arbitrary neighbourhood $X$ it is no longer true that two paths which ``cross'' have to intersect in an open point. An attempt to remedy this was made by Durrett and Schonmann \cite{Durrett87}, whose approach will be of use to us. Yet, when restricted to PPCA, their result only applies to the ones with neighbourhood of consecutive sites as above, making it trivial (their main idea is not needed for those models). As their work is somewhat informal, we indicate that this follows from the restrictive hypothesis (H3) located at the end of Sec. 4 of \cite{Durrett87}.

Improving on the approach of \cite{Durrett87} and using \cref{th:strict:monotonicity}, in \cref{app:DS} we outline how to obtain the following result.
\begin{thm}
\label{th:DS}
If for some $p\in[0,1]$ we have $\a(p)>\b(p)$, then $\theta(p)>0$ and
\begin{align*}
\lim_{p'\to p-}\a(p')={}&\a(p)&\lim_{p'\to p-}\b(p')={}&\b(p).
\end{align*}
\end{thm}
In particular, this implies $\a(\pc)\le\b(\pc)$. On the other hand, \cref{th:alpha} readily implies that $\a(p)\ge\b(p)$ for $p>\pc$. The final ingredient for proving \cref{th:2d} is the continuity to the right of $\a$ (and $\b$), which also follows from \cref{th:alpha}, since $\a$ is the decreasing limit of the continuous non-decreasing functions $\bbE_p[r_t]$.\footnote{Indeed, $\bbE_p[r_t]$ is the limit of the polynomials $\bbE_p[\max(r_t,-M)]$ as $M\to\infty$. The limit is uniform for $p\in[p',1]$ for $p'>0$, since the negative tail of $r_t$ is bounded by a geometric variable with success rate $(p')^t$. Recalling \cref{eq:subcrit} and $\pc>0$ (by comparison with branching), we further have $\a(p)=-\infty$ and $\b(p)=\infty$ for $p<\pc$.} Combining these properties, we get
\[\lim_{p\to\pc+}\a(p)-\b(p)=0,\]
which, together with \cref{th:strict:monotonicity,th:identification}, implies \cref{th:2d}.
We note that in higher dimensions it is unknown whether $\bigcap_{p>\pc}\Uc(p)$ is empty, a singleton or a larger set.

\appendix
\section{Proofs}
\label{app}
In this appendix we gather the proofs of the novel steps in the proof of the main results. The following basic result for 2dOP proved by contour arguments (see \cites{Durrett84,Durrett88}) will be used several times.
\begin{prop}
\label{prop:2dOP}
For every $\e>0$ there exist $c,\d>0$ such that for 2dOP with parameter $p\ge1-\d$ it holds that for all finite $A\subset\bbZ$ and integer $t$
\begin{align*}\bbP_p(|\{a\in A,a\to\infty\}|/|A|\le 1-\e)<{}&e^{-c|A|}\\
\bar\bbP_{p}(|\x_t^\bo|/t\le 1-\e)<{}& e^{-c t}.
\end{align*}
\end{prop}

\subsection{Primal-dual intersection---proof of Theorem~\ref{thm:GOSPproperties2}}
\label{app:coupling:preliminary}
Recall the notation of \cref{thm:GOSPproperties2,subsec:restart}.

Observe that the BG renormalisation restricts the process to a space-time slab in which all but one space dimension are suppressed. Throughout the paper we assumed this to be the~$d-1$\textsuperscript{st} dimension, but we can replace~$\be_{d-1}$ by any~$\be_i$ in \cref{thm:BGrenorm1} for $i\in\{1,\dots,d-1\}$. If the model is symmetric, the parameters $n,h,\bw,\bv$ can be chosen to be the same in all directions, however this is not necessarily the case in general, so we will need to make our notation accordingly more precise. Let us fix $\e$ and denote the parameters in \cref{thm:BGrenorm1} corresponding to~$\be_i$ by $n^{(i)},h^{(i)}\in\bbZ$, $\bw^{(i)}\in\mathbb{Z}^{d-1}$ and $\bv^{(i)}\in\mathbb{R}^{d-1}$. We then set $\bv=\sum_{i=1}^{d-1}\bv^{(i)}/(d-1)$. For simplicity we will disregard the offset $\bz_t$.

We start by noticing that by \cref{eq:property1} we may assume that $\t^A=\infty$ in $\o$ and $\tilde\t^{B+(2t\bv,0)}=\infty$ in $\o$ translated by $-2t\be_d$. We can then choose~$(\bx, s)\in A$ and~$(\tilde \bx,\tilde s)\in B+2t(\bv,1)$ such that~$\tau^{\{(\bx, s)\}}=\tilde\tau^{\{(\tilde \bx,\tilde s-2t)\}}=\infty$.

\begin{figure}
	\centering
	\begin{tikzpicture}[x=0.35cm,y=0.35cm,z=0.3cm,>=stealth]
	\draw[->] (xyz cs:x=0) -- (xyz cs:x=25);
	\draw[->] (xyz cs:y=0) -- (xyz cs:y=25);
	\draw[->] (xyz cs:z=0) -- (xyz cs:z=-13.5);
	
	\draw[thick] (-3pt,12) -- (3pt,12) node[left=6pt] {$t$};
	\draw[thick] (-3pt,24) -- (3pt,24) node[left=6pt] {$2t$};
	
	\draw[dotted] 
	(xyz cs:z=-10) -- 
	+(5,0) coordinate (u) -- 
	(xyz cs:x=5);
	\node[fill,circle,inner sep=1pt,label={below:$(\bx, s)$}] at (u) {};
	\draw[very thick] (xyz cs:x=5, z=-10) to [curve through={(xyz cs:x=5.1, y=0.2, z=-9.3)..(xyz cs:x=5.7, y=0.22, z=-9.33)}] (xyz cs:x=6, y=0.4, z=-9.5);
	\node[fill,circle,inner sep=1pt] at (xyz cs:x=6, y=0.4, z=-9.5) {};
	\draw[color=blue,thick] 
	(xyz cs:x=1.5, y=0.4, z=-10) -- 
	(xyz cs:x=11.5, y=0.4, z=-10)  coordinate (u1)  -- 
	(xyz cs:x=11.5, y=0.4, z=-9)   coordinate (v1)  -- 
	(xyz cs:x=1.5, y=0.4, z=-9)   coordinate (w1)  --
	(xyz cs:x=1.5, y=0.4, z=-10)--
	(xyz cs:x=4.5, y=5, z=-9)--
	(xyz cs:x=14.5, y=5, z=-9)   coordinate (u2)   -- 
	(xyz cs:x=14.5, y=5, z=-8)   coordinate (v2)   -- 
	(xyz cs:x=4.5, y=5, z=-8)   coordinate (w2)  --
	(xyz cs:x=4.5, y=5, z=-9);
	\draw[color=blue,thick] (u1) -- (u2);
	\draw[color=blue,thick] (v1) -- (v2);
	\draw[color=blue,thick] (w1) -- (w2);
	\foreach \n in {7,7.3,7.8, 8.1, 9, 9.4, 10, 10.4, 10.7, 11.4, 12.3, 12.6, 13.3, 13.8, 14.2, 14.8, 15.4, 15.5, 16.3}
	{
		\node[fill,circle,color=blue,inner sep=1pt] at (xyz cs:x=\n-2, y=5, z=-8.5) {};
	}
	\draw[very thick] (xyz cs:x=13.5, y=5, z=-8.5) to [curve through={(xyz cs:x=13.9, y=5.2, z=-8.3)}] (xyz cs:x=13.4, y=5.8, z=-8);
	\node[fill,circle,inner sep=1pt] at (xyz cs:x=13.4, y=5.8, z=-8) {};
	
	\draw[very thick] (xyz cs:x=14.4, y=11.8, z=-8.25) to [curve through={(xyz cs:x=14.1, y=12.3, z=-7)}] (xyz cs:x=13.8, y=12.6, z=-6.7);
	\draw[color=red,thick] 
	(xyz cs:x=12.4, y=5.8, z=-11) -- 
	(xyz cs:x=14.4, y=5.8, z=-11)  coordinate (ur1)  -- 
	(xyz cs:x=14.4, y=5.8, z=-5)   coordinate (vr1)  -- 
	(xyz cs:x=12.4, y=5.8, z=-5)   coordinate (wr1)  --
	(xyz cs:x=12.4, y=5.8, z=-11)--
	(xyz cs:x=13.4, y=11.8, z=-9)--
	(xyz cs:x=15.4, y=11.8, z=-9)   coordinate (ur2)   -- 
	(xyz cs:x=15.4, y=11.8, z=-3)   coordinate (vr2)   -- 
	(xyz cs:x=13.4, y=11.8, z=-3)   coordinate (wr2)  --
	(xyz cs:x=13.4, y=11.8, z=-9);
	\draw[color=red,thick] (ur1) -- (ur2);
	\draw[color=red,thick] (vr1) -- (vr2);
	\draw[color=red,thick] (wr1) -- (wr2);
	\foreach \n in {-3.4,-3.9,-4.2, -4.7,-5.3, -5.6,-6.1,-6.8,-7.2,-7.5,-8.25,-8.5}
	{
		\node[fill,circle,color=red,inner sep=1pt] at (xyz cs:x=14.4, y=11.8, z=\n) {};
	}
	
	\draw[dotted] (20,0) |- (0,24);
	\node[fill,circle,inner sep=1pt,label={above:$(\tilde \bx, \tilde s)$}] at (20,24) {};
	\draw[very thick] (20, 24) to [curve through={(xyz cs:x=19.8, y=23.9, z=-0.3)..(xyz cs:x=19.3, y=23.5, z=-0.7)}] (xyz cs:x=18.8, y=23.5, z=-1);
	\node[fill,circle,inner sep=1pt] at (xyz cs:x=18.8, y=23.5, z=-1) {};
	\draw[color=blue,thick] 
	(xyz cs:x=1.5+9.3, y=0.4+18.7, z=-10+7) -- 
	(xyz cs:x=11.5+9.3, y=0.4+18.7, z=-10+7)  coordinate (uu1)  -- 
	(xyz cs:x=11.5+9.3, y=0.4+18.7, z=-9+7)   coordinate (vv1)  -- 
	(xyz cs:x=1.5+9.3, y=0.4+18.7, z=-9+7)   coordinate (ww1)  --
	(xyz cs:x=1.5+9.3, y=0.4+18.7, z=-10+7)--
	(xyz cs:x=6.5+7.3, y=5+18.7, z=-9+7)--
	(xyz cs:x=16.5+7.3, y=5+18.7, z=-9+7)   coordinate (uu2)   -- 
	(xyz cs:x=16.5+7.3, y=5+18.7, z=-8+7)   coordinate (vv2)   -- 
	(xyz cs:x=6.5+7.3, y=5+18.7, z=-8+7)   coordinate (ww2)  --
	(xyz cs:x=6.5+7.3, y=5+18.7, z=-9+7);
	\draw[color=blue,thick] (uu1) -- (uu2);
	\draw[color=blue,thick] (vv1) -- (vv2);
	\draw[color=blue,thick] (ww1) -- (ww2);
	\foreach \n in {7,7.3,7.8, 8.1, 8.7, 9, 9.4, 10, 10.4, 10.7, 11.4, 12.3, 12.6, 13.3, 13.8, 14.2, 14.8, 15.4, 15.9, 16.3}
	{
		\node[fill,circle,color=blue,inner sep=1pt] at (xyz cs:x=\n+4.3, y=19.1, z=-2.5) {};
	}
	\draw[very thick] (xyz cs:x=15, y=19.1, z=-2.5) to [curve through={(xyz cs:x=15.1, y=18.9, z=-2.8)}] (xyz cs:x=14.8, y=18.6, z=-3);
	\node[fill,circle,inner sep=1pt] at (xyz cs:x=14.8, y=18.6, z=-3) {};
	\draw[color=red,thick] 
	(xyz cs:x=12.4+0.4, y=5.8+6.8, z=-11+3) -- 
	(xyz cs:x=14.4+0.4, y=5.8+6.8, z=-11+3)  coordinate (uur1)  -- 
	(xyz cs:x=14.4+0.4, y=5.8+6.8, z=-5+3)   coordinate (vvr1)  -- 
	(xyz cs:x=12.4+0.4, y=5.8+6.8, z=-5+3)   coordinate (wwr1)  --
	(xyz cs:x=12.4+0.4, y=5.8+6.8, z=-11+3)--
	(xyz cs:x=13.4+0.4, y=11.8+6.8, z=-9+3)--
	(xyz cs:x=15.4+0.4, y=11.8+6.8, z=-9+3)   coordinate (uur2)   -- 
	(xyz cs:x=15.4+0.4, y=11.8+6.8, z=-3+3)   coordinate (vvr2)   -- 
	(xyz cs:x=13.4+0.4, y=11.8+6.8, z=-3+3)   coordinate (wwr2)  --
	(xyz cs:x=13.4+0.4, y=11.8+6.8, z=-9+3);
	\draw[color=red,thick] (uur1) -- (uur2);
	\draw[color=red,thick] (vvr1) -- (vvr2);
	\draw[color=red,thick] (wwr1) -- (wwr2);
	\foreach \n in {-3.3,-3.6,-4.1, -4.3, -4.9,-5.3, -5.7,-6.1,-6.5,-7.27,-7.7,-8.1,-8.45}
	{
		\node[fill,circle,color=red,inner sep=1pt] at (xyz cs:x=13.8, y=12.6, z=\n+1) {};
	}
	\end{tikzpicture}
		\caption{The argument for $d=3$. The bottom and top tilted boxes correspond to the $\be_2$ space-time slabs, while the middle two correspond to $\be_1$ space-time slabs. Dots represent
		infected translates of the boxes $B_{n^{(2)}}$ and $B_{n^{(1)}}$ respectively.}
	\label{fig:primaldual}
\end{figure}
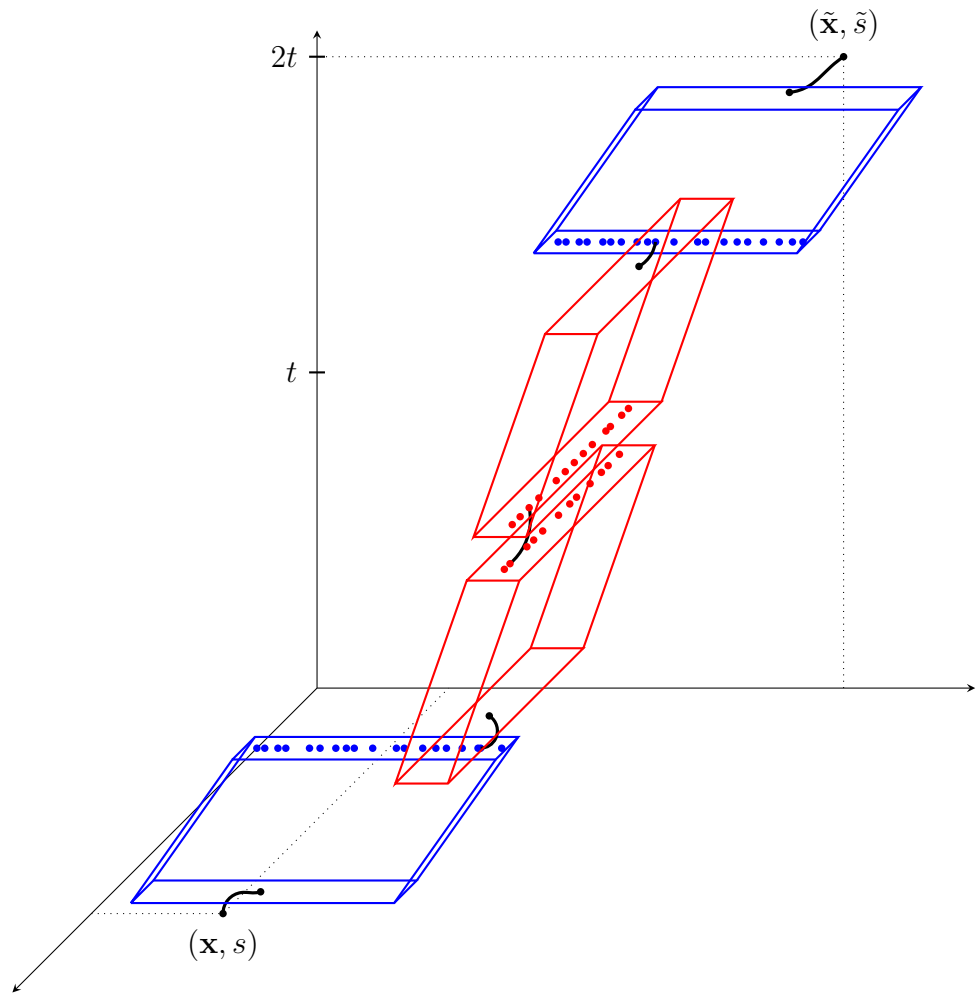

We can perform a restart argument starting from $(\bx,s)$ and $(\tilde \bx,\tilde s)$ until they simultaneously infect a (translate of the) box $B_{n^{(d-1)}}$ each and that the two boxes give rise to a percolating 2dOP in their respective $\be_{d-1}$-space-time slabs. As in \cref{subsec:restart} the restart argument is exponentially unlikely to require more than $\d t$ steps for some small $\d>0$, the positions of the two boxes differ by $2t(\bv,1)+\bD$ with $\|\bD\|=O(c+\d)t$ (here asymptotic notation is w.r.t.\ $t\to\infty$).

Informally, the rest of the argument is as follows (see \cref{fig:primaldual}). We let the primal and dual processes evolve in their respective $\be_{d-1}$-space-time slabs for $t/(d-1)$ time steps. The comparison with the percolating 2dOP and \cref{prop:2dOP} then ensure that both processes reach many infected copies of the box $B_{n^{(d-1)}}$ at times $t/(d-1)$ and $2t-t/(d-1)$ respectively. As the space-time slab is in the $\be_{d-1}$ direction, we can find a lot of ``well-aligned'' pairs of infected primal and dual boxes. Namely, we require their $d-1$\textsuperscript{st} coordinates to differ exactly by the amount of tilting we will have in the rest of the procedure: $\frac{2t}{d-1}\sum_{i=1}^{d-2} v^{(i)}_{d-1}$. Roughly speaking, at this point we have managed to cancel the $d-1$\textsuperscript{st} coordinate of $\bD$. With high probability, at least one such pair infects a (translate of the) box $B_{n^{(d-2)}}$ each, that give rise to a percolating 2dOP in their respective $\be_{(d-2)}$-space-time slabs. We then repeat the same reasoning for each direction. Eventually, at time $t$, we find many infected pairs of primal and dual boxes sufficiently close to each other, so that with high probability there will be an open path between at least one pair.

More precisely, fix a large integer $K$ so that \cref{eq:l} holds for $n=n^{(d-2)}$ and $t=\lfloor Kh^{(d-1)}/2\rfloor$. Then by \cref{prop:2dOP} we get that at time $t/(7h^{(d-1)}(d-1))-K$ both in $\z$ and $\tilde\z$ (the renormalised 2dOP corresponding to $\x$ and $\tilde\x$) infect at least $2/3$ of the (renormalised) sites that can be reached from the sites corresponding to the initial two boxes $B_{n^{(d-1)}}$ close to $(\bx,s)$ and $(\tilde\bx,\tilde s)$. By the pigeonhole principle, as $c$ and $\d$ are sufficiently small, there are at least $t/(22h^{(d-1)}(d-1))$ sites $(z,t/(7h^{(d-1)}(d-1))-K)$ which are infected in $\z$ and such that $(\tilde z,t/(7h^{(d-1)}(d-1))-K)$ is infected in $\tilde\z$ with $\tilde z-z=-\lfloor\D_{d-1}/w^{(d-1)}_{d-1}\rfloor$. It then follows from~\cref{eq:l}, \cref{prop:2dOP} and the pigeonhole principle that up to an exponentially unlikely event at least one such couple $z,\tilde z$ gives rise to two boxes $B_{n^{(d-2)}}+(\by,t/(d-1))$ and $B_{n^{(d-2)}}+(\tilde \by,t(2-1/(d-1)))$ infected in $\x$ and $\tilde \x$ respectively such that $\tilde \by-\by=2tv-2t\bv^{(d-1)}/(d-1)+\sum_{i=1}^{d-2}\D_i\be_i$ and such that the 2dOP renormalisations in direction $\be_{d-2}$ of each of the boxes percolate.

Repeating the same reasoning for each direction and recalling the definition of $\bv$, we obtain the desired conclusion.

\subsection{Tilting}
\label{app:tilting}
Recall the setting of \cref{subsec:shape}. In this section we show how to remove the additional assumption $\bv=0$ used there in the proofs of \cref{eq:main:lower,eq:main:upper} and the continuity of $U$ in \cref{th:main}, as well as \cref{thm:largedev}. The reasoning for \cref{thm:largedev} and the continuity being identical to the one for \cref{eq:main:lower,eq:main:upper}, we only address the latter.

Indeed, we can assume w.l.o.g.\ that the vector $\bv$ in \cref{thm:GOSPproperties3} is in $\bbQ^{d-1}$ and then apply the linear map $(\bx,t)\mapsto(\bx-t\bv,t)$ to the lattice. We will refer to the resulting lattice $\hat\bbZ^d$ as the \emph{tilted lattice} and define its \emph{period} $\hat R:=\min\{t\in\bbZ:t\bv\in\bbZ^{d-1},t\ge R\}$ and \emph{base} $\hat B=(\bbR^{d-1}\times[0,\hat R))\cap\hat\bbZ^d$. For $A\subset S$ we define the tilted process, hitting time, hit and coupled regions
\begin{align*}
\hat\x_t^A:={}&\left\{(\bx,s)\in\hat B:(\bx+(t+s)\bv,0)\in\x_{t+s}^A\right\},\\
\hat K_t^A:={}&\left\{(\bx,s)\in \hat B:\hat\x^A_t(\bx,s)=\hat\x^S_t(\bx,s)\right\},\\
\hat t^A(\bx,s):={}&\min\left\{t\in s+\hat R\bbZ:t\ge0,(\bx+t\bv,0)\in\x^A_t\right\}\quad (\bx,s)\in\hat B,\\
\hat H_t^A:={}&\left\{(\bx,s)\in \hat B: \hat t^A(\bx,s)\le t\right\}
\end{align*}

\begin{figure}
	\begin{center}
	\begin{subfigure}{0.45\textwidth}
	\centering
	\begin{tikzpicture}[line cap=round,line join=round,>=triangle 45,x=0.1\textwidth,y=0.05\textwidth]
		\draw[->] (-9,0) -- (1,0);
		\draw[->] (-8,-1) -- (-8,10);
		\draw (-8,0) -- (-6,9);
		\draw (-8,0) -- (-1,9);
		\draw (-4.3,8) -- (-7.855,0);
		\fill[fill=black,fill opacity=0.3] (-8,0) -- (-4.8,9) -- (-3.22,9) -- cycle;
		\fill [color=black] (-4.3,8) circle (1pt);
		\fill [color=black] (-7.855,0) circle (1pt);

		
		\begin{scriptsize}
		\draw (-4.3,8) node [above right] {$(x,t)$};
		\draw (-7.855,0)--(-7.2,0.4) node [right] {$(x-vt,0)$};
		\draw (-4.89,4)--(-4.89,3.6) node [anchor=north west,xshift=-0.1\textwidth] {$\bigcup_{t>0}\left((1-2\varepsilon)t\hat U+vt\right)\times\{t\}$};
		\draw (-8.5,-0.1)--(-8.5,0.1) node [below=0.1cm] {$-ct$};
		\draw (-7.5,-0.1)--(-7.5,0.1) node [below=0.1cm] {$ct$};
		\end{scriptsize}
		\end{tikzpicture}
		\caption{\label{fig:shapethm1}Case $\|x-vt\|\le ct$.}
	\end{subfigure}
	\quad
	\begin{subfigure}{0.45\textwidth}
	\centering
	\begin{tikzpicture}[line cap=round,line join=round,>=triangle 45,x=0.1\textwidth,y=0.05\textwidth]
		\draw[->] (3,0) -- (13,0);
		\draw[->] (4,-1) -- (4,10);
	    \draw (4,0) -- (6,9);
		\draw (4,0) -- (11,9);
		\draw (10,8.5) -- (6.3333,0);
		\fill[fill=black,fill opacity=0.3] (7,3.05) -- (9.42,9) -- (10.57,9) -- cycle;
		\fill [color=black] (10,8.5) circle (1pt);
		\fill [color=black] (6.3333,0) circle (1pt);
		\fill [color=black] (8.28,4.5) circle (1pt);
		\fill[pattern=north west lines] (4,0) -- (12, 9) -- (13, 9) -- cycle;
		\fill[pattern=north west lines] (4,0) -- (5.5, 9) -- (4.9, 9) -- cycle;
		\draw[color=red,very thick] (4,0) to [curve through={(4.3,1.1)..(7,3.05)..(7.9,4.1)}] (8.28,4.5);
		\fill [color=black] (8.28,4.5) circle (1pt);
		\fill [color=black] (7,3.05) circle (1pt);
		\begin{scriptsize}
		\draw (10,8.5) node [left] {$(x,t)$};
		\draw (8.28,4.5) node [right=0.1cm] {$(y,s)$};
		\draw (7,3.05) node [below right] {$(y',s')$};
		\draw (6.3333,0) node [below] {$(x-vt,0)$};
		\end{scriptsize}
		\end{tikzpicture}
		\caption{\label{fig:shapethm2}Case $\|x-vt\|> ct$.}
		\end{subfigure}
		\caption{The original lattice for~$d=2$. Shaded areas represent the cone $\bigcup_{t>0} [t(v-c),t(v+c)]\times\{t\}$ rooted at $o$ and~$(y', s')$ respectively. The hatched region is $\D$.} 
		\label{fig:shapethm}
	\end{center}
\end{figure}
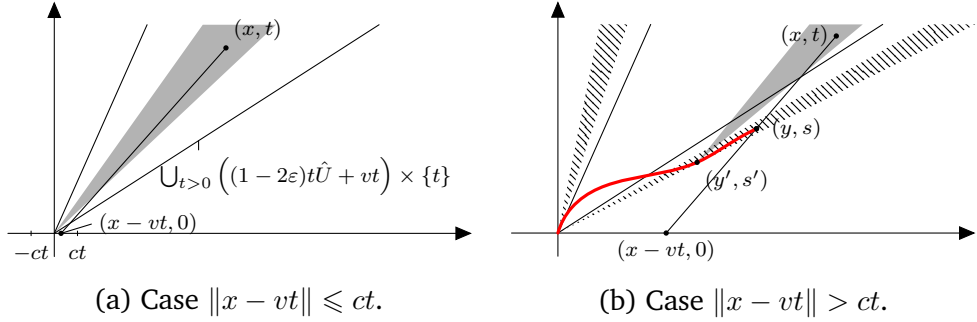

The proof of \cref{subsec:shape} applies to GOSP in the tilted setting, yielding \cref{eq:main:lower,eq:main:upper} in $\hat\bbZ^d$ for some convex compact limit shape $\hat U\subset \bbR^{d-1}$ containing $o$ in its interior. We then need to transfer the result back to the original lattice with $U=\hat U+\bv$. By the definition of~$\hat\x_t$ and~$\hat K_t$, \cref{eq:main:upper} and the inclusion in the coupled region in~\cref{eq:main:lower} are immediate. It remains to show that for every~$\e>0$,~$\bar\bbP_p$-a.s.\ for every~$t$ large enough~$H_t{}\supset(((1-\e)tU)\times[0,R))\cap\bbZ^d$. Our strategy, somewhat similar to \cite{Durrett82}, is as follows. Fixing $\bx\in\bbZ^{d-1}$ such that $(\bx,0)$ should belong to $H_t$, we trace the line of slope $\bv$ from $(\bx,t)$ and determine when it intersects the boundary of the cone $\bigcup_{t'\ge 0} (t'U)\times\{t'\}$ (see \cref{fig:shapethm}). Someone close to the intersection point should be infected around time $t'$ by the result available in $\hat \bbZ^d$. But then at a time $t'-\e t$ some site close to the intersection has survived for time $\e t$, so, applying \cref{eq:property1,thm:GOSPproperties3} to this site we manage to reach $(\bx,t)$ as desired. With this in mind, let us spell out the details.

\Cref{eq:main:lower,eq:main:upper} of \cref{th:main} on $\hat\bbZ^d$ imply that for every~$\e>0$ $\bar\bbP_p$-a.s.\ there exists a constant~$C$ such that for every~$(\bx, s)\in\hat B$ with~$\| \bx \|\geq C$
\begin{equation}\label{eq:boundary}
\left(\bx, \hat t(\bx, s)\right)\in \bigcup_{t>0}\left(t\partial \hat U\right)\times\{t\},
\end{equation}
setting $\partial\hat U:=((1+\e)\hat U)\setminus ((1-\e)\hat U)$. Observe that this event implies that in the original lattice there is at least one vertex infected by the origin in the intersection of $\D:=\bigcup_{t>0}(t\partial\hat U_1+t\bv)\times\{t\}$ and the ray $(\bx+t\bv,t)_{t\ge s}$ for all $(\bx,s)\in\hat B$ such that $\|\bx\|\ge C$ (see \cref{fig:shapethm2}).

Fix~$c$ and~$\e$ so that \cref{thm:GOSPproperties,thm:GOSPproperties3} hold for the original lattice. We now argue that for any~$t>C/c$ and for any~$\bx\in ((1-2\e)t\hat U+vt)\cap\mathbb Z^{d-1}$, we have~$(\bx, 0)\in H_t$ except with probability exponentially small in~$t$. If~$\|\bx-vt\|\le ct$ (see \cref{fig:shapethm1}), then~\cref{thm:GOSPproperties3} directly gives the desired result.

Assuming~$\|\bx-t\bv\|>ct$, let~$\d>0$ be small enough depending on $X,c,\hat U,\bv,\e$, but not $t,\bx,C$. \Cref{eq:boundary} implies that~$\bar\bbP_p$-a.s.\ there is a vertex~$(\by, s)$ in the intersection of~$\D$ and the segment from~$(\bx-t\bv,0)$ to~$(\bx,t)$, such that~$(\by,0)\in\x_s^\bo$. As $s
>\delta t$, we can take a site along an associated infection path at time closest to~$s-\delta t$, and denote it by~$(\by', s')$ (see \cref{fig:shapethm2}).

We then have that~$(\by',0)\in\x^\bo_{s'}$ and~$(\by', s')$ survives for time at least~$\delta t/2$, so we use~\cref{eq:property1} to conclude that~$(\by', s')$ survives with probability exponentially high in $t$. Once we have survival, we can use~\cref{thm:GOSPproperties3} to show that~$\bx$ is in the hit region of~$(\by', s')$ at time~$t$ with high probability, thus it is also in the hit region of the origin. Indeed, \[\left\|\bx-\by'-(t-s')\bv\right\|=\left\|\by-\by'-(s-s')\bv\right\|\le t\sqrt{\d}<c(t-s'),\]
since $(\bx,t)$ is at distance at least $\k t$ from $\D$ (and thus from $(\by,s)$) for some $\k>0$ depending only on $\hat U,\bv,\e$.

This completes the proof of~\cref{eq:main:lower,eq:main:upper} of~\cref{th:main} as stated for the original lattice.

\subsection{Density large deviations---proof of Eq.~(\ref{eq:LD:reduced})}
\label{app:LD}
Recall the notation of \cref{thm:LD} and \cref{subsec:cluster:properties}. We will use a similar argument to \cite{Durrett88} but based on a completely different renormalisation. 

Let us fix $p>\pc$ and $a<\tilde\theta(p)$, let $C$ be large enough depending on $p$, let $L$ be large enough depending on $p,a,C$ and define $\bw=(L,\dots,L)\in\bbZ^{d-1}$ and $s=CL+L/C$. Recalling \cref{eq:Bwhv}, let $B=B(\bw,CL,\bv)$ with $\bv$ as in \cref{thm:GOSPproperties3}. We say that $B$ is \emph{good} if the following events all occur.
\begin{enumerate}
\item\label{event1} For each site $(\bx,t)\in B\cap S=B(\bw,R,\bv)$ we have either $\t^{(\bx,t)}<L/C$ or $\t^{(\bx,t)}\ge s$, where $\t^{(\bx,t)}$ is defined as $\t^{\{(\bx,R-1)\}}$ for the configuration $\o$ translated by $-(t-R+1)\be_d$.
\item\label{event2} For each site $(\bx,t)\in B\cap S$ such that $\t^{(\bx,t)}\ge L/C$ we have $K_{s}^{(\bx,t)}\supset B(3\bw,R,\bv)+s\bv$ and $K^{(\bx,t)}_{CL}\supset B(3\bw,R,\bv)+CL\bv$ with $K_{u}^{(\bx,t)}$ defined as $K^{\{(\bx,R-1)\}}_{u-t+R-1}$ for the configuration $\o$ translated by $-(t-R+1)\be_d$.
\item \label{event3} \makebox[\linegoal]{$\displaystyle\begin{aligned}[t]
\x^S_{s}\cap \left(B(\bw/C,R,\bv)+s(v,0)+L\be_{d-1}\right)&{}\neq \varnothing,\\
\x^S_s\cap \left(B(\bw/C,R,\bv)+s(\bv,0)-L\be_{d-1}\right)&{}\neq \varnothing.
\end{aligned}$}
\end{enumerate}
In words, each site which does not die quickly survives well beyond the top of $B$ and infects the same set of sites at the top of $B\pm L\be_{d-1}$, at least one of which does not die quickly. Indeed, the neighbourhood $X$ being finite, the only way to reach $B(\bw/C,R,\bv)+s(\bv,1)\pm L\be_{d-1}$ is to go through $B(\bw,R,\bv)+CL(\bv,1)\pm L\be_{d-1}$. Therefore, considering a renormalised two-dimensional lattice with sites corresponding to disjoint translates of $B$, the resulting 2dOP is $C^{2}$-dependent, as $B$ being good only depends on the configuration in $B(C^2\bw,2CL,\bv)$.

We next show that the parameter of the 2dOP is close to $1$ when $L$ is large enough, so that by \cite{Liggett97} it stochastically dominates an independent 2dOP with parameter close to 1. Indeed, Event \ref{event1} fails with exponentially small (in $L$) probability by \cref{eq:property1}; Event \ref{event2} fails with exponentially small probability by \cref{thm:GOSPproperties3,eq:property1}; Event \ref{event3} fails with stretched exponentially small probability by \cref{eq:property2} applied to the dual process.

It is easily checked that if a renormalised site $B$ percolates in 2dOP, then each site in $B\cap S$ either dies in time at most $L/C$ or also percolates. Recalling \cref{prop:2dOP}, the rest of the proof is essentially as in \cite{Durrett88}. Taking $n$ much larger than $L$, we may cut a box $B_n$ into $(n/(2L))^{d-2}$ strips each giving rise to a different renormalised 2dOP. It is then standard to show that the total proportion of percolating renormalised sites is not close to 1 with probability at most $\exp(-\e n^{d-1})$ for some $\e>0$ depending on $L$ but not on $n$. Moreover, by standard large deviations for independent random variables, the proportion of sites, which survive at least $L/C$ steps in $B_n$ is smaller than $\theta(p)$ with probability at most $\exp(-\e'n^{d-1})$ for some $\e'>0$ depending on $L$ but not on $n$. We may then conclude by discarding the renormalised sites which do not percolate. Finally, performing the same reasoning for the dual process rather than the primal one, we obtain the desired conclusion (with $\tilde\theta(p)$ instead of  $\theta(p)$).

\subsection{Enhanced 2d renormalisation---proof of Theorem~\ref{th:DS}}
\label{app:DS}

Recall the notation of \cref{th:DS} and \cref{sec:2d}.

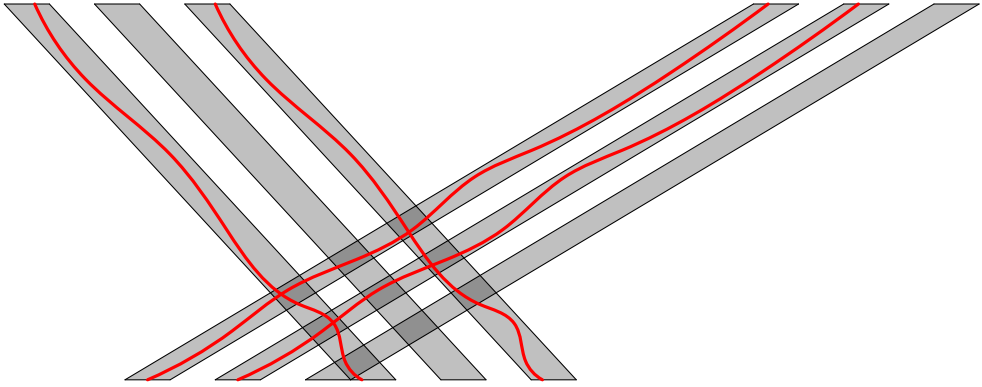
\begin{figure}
    \centering
\begin{tikzpicture}[line cap=round,line join=round,>=triangle 45,x=0.75cm,y=0.5cm]
\fill[fill=black,fill opacity=0.25] (-0.4,0) -- (0.4,0) -- (11.54,10) -- (10.74,10) -- cycle;
\fill[fill=black,fill opacity=0.25] (1.2,0) -- (2,0) -- (13.14,10) -- (12.34,10) -- cycle;
\fill[fill=black,fill opacity=0.25] (2.8,0) -- (3.6,0) -- (14.74,10) -- (13.94,10) -- cycle;
\fill[fill=black,fill opacity=0.25] (-2.54,10) -- (-1.74,10) -- (4.4,0) -- (3.6,0) -- cycle;
\fill[fill=black,fill opacity=0.25] (-0.94,10) -- (-0.14,10) -- (6,0) -- (5.2,0) -- cycle;
\fill[fill=black,fill opacity=0.25] (0.66,10) -- (1.46,10) -- (7.6,0) -- (6.8,0) -- cycle;
\draw (-0.4,0)-- (0.4,0);
\draw (0.4,0)-- (11.54,10);
\draw (11.54,10)-- (10.74,10);
\draw (10.74,10)-- (-0.4,0);
\draw (1.2,0)-- (2,0);
\draw (2,0)-- (13.14,10);
\draw (13.14,10)-- (12.34,10);
\draw (12.34,10)-- (1.2,0);
\draw (2.8,0)-- (3.6,0);
\draw (3.6,0)-- (14.74,10);
\draw (14.74,10)-- (13.94,10);
\draw (13.94,10)-- (2.8,0);
\draw (-2.54,10)-- (-1.74,10);
\draw (-1.74,10)-- (4.4,0);
\draw (4.4,0)-- (3.6,0);
\draw (3.6,0)-- (-2.54,10);
\draw (-0.94,10)-- (-0.14,10);
\draw (-0.14,10)-- (6,0);
\draw (6,0)-- (5.2,0);
\draw (5.2,0)-- (-0.94,10);
\draw (0.66,10)-- (1.46,10);
\draw (1.46,10)-- (7.6,0);
\draw (7.6,0)-- (6.8,0);
\draw (6.8,0)-- (0.66,10);
\draw[color=red,very thick] (0,0) to [curve through={(1.5,1.3)..(2.5,2.4)..(3.7,3.2)..(4.7,4)..(5.8,5.4)..(6.8,6.1)..(9.3,8.2)}] (11,10);
\draw[color=red,very thick] (1.6,0) to [curve through={(3.1,1.3)..(4.1,2.4)..(5.3,3.2)..(6.3,4)..(7.4,5.4)..(8.4,6.1)..(10.9,8.2)}] (12.6,10);
\draw[color=red,very thick] (3.8,0) to [curve through={(3.5,0.5)..(3.3,1.5)..(2.7,2)..(1.4,4)..(0.4,6)..(-1.1,8)}] (-2,10);
\draw[color=red,very thick] (7,0) to [curve through={(6.7,0.5)..(6.5,1.5)..(5.9,2)..(4.6,4)..(3.6,6)..(2.1,8)}] (1.2,10);
\end{tikzpicture}
    \caption{The shaded boxes are likely to contain paths crossing them. In order to transition from one path to another, we use additional infections as illustrated in \cref{fig:DS2}.}
    \label{fig:DS1}
\end{figure}
The idea of \cite{Durrett87} is to introduce several translates of the original long boxes from \cref{eq:crossing} and hope that their paths have positive probability of intersecting and do so independently (see \cref{fig:DS1}). More precisely, increasing $p$ by a small amount $\epsilon$, they require that the additional infections suffice to transition with positive probability from one path to the other at the place where they cross. As we shall see, although it is not possible to do this, as intended, in one step at the crossing point, it is possible to find a place where to do it in several steps.

\begin{figure}
    \centering
\begin{tikzpicture}[line cap=round,line join=round,>=triangle 45,x=1.0cm,y=0.65cm]

\fill[fill=black,fill opacity=0.25] (-1.76,1.22) -- (-1.76,0.72) -- (-0.96,0.72) -- (-0.96,1.22) -- cycle;
\fill[fill=black,fill opacity=0.25] (-2.08,1.6) -- (-2.08,1.1) -- (-1.28,1.1) -- (-1.28,1.6) -- cycle;
\fill[fill=black,fill opacity=0.25] (-1.85,1.93) -- (-1.85,1.43) -- (-1.05,1.43) -- (-1.05,1.93) -- cycle;
\fill[fill=black,fill opacity=0.25] (-1.26,2.2) -- (-1.26,1.7) -- (-0.46,1.7) -- (-0.46,2.2) -- cycle;
\fill[fill=black,fill opacity=0.25] (-0.7,2.65) -- (-0.7,2.15) -- (0.1,2.15) -- (0.1,2.65) -- cycle;
\fill[fill=black,fill opacity=0.25] (-0.21,2.93) -- (-0.21,2.43) -- (0.59,2.43) -- (0.59,2.93) -- cycle;
\fill[fill=black,fill opacity=0.25] (0.48,3.7) -- (0.48,3.2) -- (1.28,3.2) -- (1.28,3.7) -- cycle;
\fill[fill=black,fill opacity=0.25] (0.38,3.23) -- (0.38,2.73) -- (1.18,2.73) -- (1.18,3.23) -- cycle;
\fill[fill=black,fill opacity=0.25] (0.74,4.2) -- (0.74,3.7) -- (1.54,3.7) -- (1.54,4.2) -- cycle;
\fill[fill=black,fill opacity=0.25] (1.34,4.48) -- (1.34,3.98) -- (2.14,3.98) -- (2.14,4.48) -- cycle;
\fill[fill=black,fill opacity=0.25] (2.28,5.49) -- (2.28,4.99) -- (3.08,4.99) -- (3.08,5.49) -- cycle;
\fill[fill=black,fill opacity=0.25] (1.89,4.86) -- (1.89,4.36) -- (2.69,4.36) -- (2.69,4.86) -- cycle;
\fill[fill=black,fill opacity=0.25] (2.91,5.9) -- (2.91,5.4) -- (3.71,5.4) -- (3.71,5.9) -- cycle;
\fill[fill=black,fill opacity=0.25] (3.27,6.3) -- (3.27,5.8) -- (4.07,5.8) -- (4.07,6.3) -- cycle;
\fill[fill=black,fill opacity=0.25] (4.04,6.73) -- (4.04,6.23) -- (4.84,6.23) -- (4.84,6.73) -- cycle;
\fill[fill=black,fill opacity=0.25] (4.29,7.17) -- (4.29,6.67) -- (5.09,6.67) -- (5.09,7.17) -- cycle;
\fill[fill=black,fill opacity=0.25] (5.21,8.07) -- (5.21,7.57) -- (6.01,7.57) -- (6.01,8.07) -- cycle;
\fill[fill=black,fill opacity=0.25] (4.45,7.65) -- (4.45,7.15) -- (5.25,7.15) -- (5.25,7.65) -- cycle;
\fill[fill=black,fill opacity=0.25] (6.76,9.46) -- (6.76,8.96) -- (7.56,8.96) -- (7.56,9.46) -- cycle;
\fill[fill=black,fill opacity=0.25] (5.71,8.56) -- (5.71,8.06) -- (6.51,8.06) -- (6.51,8.56) -- cycle;
\fill[fill=black,fill opacity=0.25] (6.51,8.96) -- (6.51,8.46) -- (7.31,8.46) -- (7.31,8.96) -- cycle;
\fill[fill=black,fill opacity=0.25] (7.47,9.67) -- (7.47,9.17) -- (8.27,9.17) -- (8.27,9.67) -- cycle;
\fill[fill=black,fill opacity=0.25] (8.09,9.97) -- (8.09,9.47) -- (8.89,9.47) -- (8.89,9.97) -- cycle;
\fill[fill=black,fill opacity=0.25] (8.37,10.43) -- (8.37,9.93) -- (9.17,9.93) -- (9.17,10.43) -- cycle;
\fill[fill=black,fill opacity=0.25] (8.86,10.71) -- (8.86,10.21) -- (9.66,10.21) -- (9.66,10.71) -- cycle;
\fill[fill=black,fill opacity=0.25] (1.71,5.12) -- (1.71,4.62) -- (2.51,4.62) -- (2.51,5.12) -- cycle;
\draw (-3.09,10)-- (-1.09,10);
\draw (-1.09,10)-- (3,0);
\draw (3,0)-- (1,0);
\draw (1,0)-- (-3.09,10);
\draw (-3,0)-- (-1,0);
\draw (-1,0)-- (10.47,10);
\draw (10.47,10)-- (8.47,10);
\draw (8.47,10)-- (-3,0);
\draw (-1.5,0.02)-- (-1.81,0.4);
\draw (-1.81,0.4)-- (-1.59,0.73);
\draw (-1.59,0.73)-- (-1,1);
\draw (-1,1)-- (-0.44,1.45);
\draw (-0.44,1.45)-- (0.05,1.73);
\draw (0.05,1.73)-- (0.65,2.03);
\draw (0.65,2.03)-- (0.74,2.5);
\draw (0.74,2.5)-- (1,3);
\draw (1,3)-- (1.6,3.28);
\draw (1.6,3.28)-- (2.15,3.66);
\draw (2.15,3.66)-- (1.97,3.92);
\draw (1.97,3.92)-- (2.54,4.29);
\draw (2.54,4.29)-- (3.18,4.7);
\draw (3.18,4.7)-- (3.53,5.1);
\draw (3.53,5.1)-- (4.3,5.53);
\draw (4.3,5.53)-- (4.55,5.97);
\draw (4.55,5.97)-- (4.71,6.45);
\draw (4.71,6.45)-- (5.47,6.87);
\draw (5.47,6.87)-- (5.98,7.36);
\draw (5.98,7.36)-- (6.78,7.76);
\draw (6.78,7.76)-- (7.03,8.26);
\draw (7.03,8.26)-- (7.73,8.47);
\draw (7.73,8.47)-- (8.35,8.77);
\draw (8.35,8.77)-- (8.64,9.23);
\draw (8.64,9.23)-- (9.12,9.51);
\draw (9.12,9.51)-- (9.48,9.99);
\draw (-1.76,1.22)-- (-1.76,0.72);
\draw (-1.76,0.72)-- (-0.96,0.72);
\draw (-0.96,0.72)-- (-0.96,1.22);
\draw (-0.96,1.22)-- (-1.76,1.22);
\draw (-2.08,1.6)-- (-2.08,1.1);
\draw (-2.08,1.1)-- (-1.28,1.1);
\draw (-1.28,1.1)-- (-1.28,1.6);
\draw (-1.28,1.6)-- (-2.08,1.6);
\draw (-1.85,1.93)-- (-1.85,1.43);
\draw (-1.85,1.43)-- (-1.05,1.43);
\draw (-1.05,1.43)-- (-1.05,1.93);
\draw (-1.05,1.93)-- (-1.85,1.93);
\draw (-1.26,2.2)-- (-1.26,1.7);
\draw (-1.26,1.7)-- (-0.46,1.7);
\draw (-0.46,1.7)-- (-0.46,2.2);
\draw (-0.46,2.2)-- (-1.26,2.2);
\draw (-0.7,2.65)-- (-0.7,2.15);
\draw (-0.7,2.15)-- (0.1,2.15);
\draw (0.1,2.15)-- (0.1,2.65);
\draw (0.1,2.65)-- (-0.7,2.65);
\draw (-0.21,2.93)-- (-0.21,2.43);
\draw (-0.21,2.43)-- (0.59,2.43);
\draw (0.59,2.43)-- (0.59,2.93);
\draw (0.59,2.93)-- (-0.21,2.93);
\draw (0.48,3.7)-- (0.48,3.2);
\draw (0.48,3.2)-- (1.28,3.2);
\draw (1.28,3.2)-- (1.28,3.7);
\draw (1.28,3.7)-- (0.48,3.7);
\draw (0.38,3.23)-- (0.38,2.73);
\draw (0.38,2.73)-- (1.18,2.73);
\draw (1.18,2.73)-- (1.18,3.23);
\draw (1.18,3.23)-- (0.38,3.23);
\draw (0.74,4.2)-- (0.74,3.7);
\draw (0.74,3.7)-- (1.54,3.7);
\draw (1.54,3.7)-- (1.54,4.2);
\draw (1.54,4.2)-- (0.74,4.2);
\draw (1.34,4.48)-- (1.34,3.98);
\draw (1.34,3.98)-- (2.14,3.98);
\draw (2.14,3.98)-- (2.14,4.48);
\draw (2.14,4.48)-- (1.34,4.48);
\draw (2.28,5.49)-- (2.28,4.99);
\draw (2.28,4.99)-- (3.08,4.99);
\draw (3.08,4.99)-- (3.08,5.49);
\draw (3.08,5.49)-- (2.28,5.49);
\draw (1.89,4.86)-- (1.89,4.36);
\draw (1.89,4.36)-- (2.69,4.36);
\draw (2.69,4.36)-- (2.69,4.86);
\draw (2.69,4.86)-- (1.89,4.86);
\draw (2.91,5.9)-- (2.91,5.4);
\draw (2.91,5.4)-- (3.71,5.4);
\draw (3.71,5.4)-- (3.71,5.9);
\draw (3.71,5.9)-- (2.91,5.9);
\draw (3.27,6.3)-- (3.27,5.8);
\draw (3.27,5.8)-- (4.07,5.8);
\draw (4.07,5.8)-- (4.07,6.3);
\draw (4.07,6.3)-- (3.27,6.3);
\draw (4.04,6.73)-- (4.04,6.23);
\draw (4.04,6.23)-- (4.84,6.23);
\draw (4.84,6.23)-- (4.84,6.73);
\draw (4.84,6.73)-- (4.04,6.73);
\draw (4.29,7.17)-- (4.29,6.67);
\draw (4.29,6.67)-- (5.09,6.67);
\draw (5.09,6.67)-- (5.09,7.17);
\draw (5.09,7.17)-- (4.29,7.17);
\draw (5.21,8.07)-- (5.21,7.57);
\draw (5.21,7.57)-- (6.01,7.57);
\draw (6.01,7.57)-- (6.01,8.07);
\draw (6.01,8.07)-- (5.21,8.07);
\draw (4.45,7.65)-- (4.45,7.15);
\draw (4.45,7.15)-- (5.25,7.15);
\draw (5.25,7.15)-- (5.25,7.65);
\draw (5.25,7.65)-- (4.45,7.65);
\draw (6.76,9.46)-- (6.76,8.96);
\draw (6.76,8.96)-- (7.56,8.96);
\draw (7.56,8.96)-- (7.56,9.46);
\draw (7.56,9.46)-- (6.76,9.46);
\draw (5.71,8.56)-- (5.71,8.06);
\draw (5.71,8.06)-- (6.51,8.06);
\draw (6.51,8.06)-- (6.51,8.56);
\draw (6.51,8.56)-- (5.71,8.56);
\draw (6.51,8.96)-- (6.51,8.46);
\draw (6.51,8.46)-- (7.31,8.46);
\draw (7.31,8.46)-- (7.31,8.96);
\draw (7.31,8.96)-- (6.51,8.96);
\draw (7.47,9.67)-- (7.47,9.17);
\draw (7.47,9.17)-- (8.27,9.17);
\draw (8.27,9.17)-- (8.27,9.67);
\draw (8.27,9.67)-- (7.47,9.67);
\draw (8.09,9.97)-- (8.09,9.47);
\draw (8.09,9.47)-- (8.89,9.47);
\draw (8.89,9.47)-- (8.89,9.97);
\draw (8.89,9.97)-- (8.09,9.97);
\draw (8.37,10.43)-- (8.37,9.93);
\draw (8.37,9.93)-- (9.17,9.93);
\draw (9.17,9.93)-- (9.17,10.43);
\draw (9.17,10.43)-- (8.37,10.43);
\draw (8.86,10.71)-- (8.86,10.21);
\draw (8.86,10.21)-- (9.66,10.21);
\draw (9.66,10.21)-- (9.66,10.71);
\draw (9.66,10.71)-- (8.86,10.71);
\draw (1.71,5.12)-- (1.71,4.62);
\draw (1.71,4.62)-- (2.51,4.62);
\draw (2.51,4.62)-- (2.51,5.12);
\draw (2.51,5.12)-- (1.71,5.12);
\draw (1.86,0)-- (2,0.17);
\draw (2,0.17)-- (1.77,0.4);
\draw (1.77,0.4)-- (1.43,0.85);
\draw (1.43,0.85)-- (1.64,1.01);
\draw (1.64,1.01)-- (1.4,1.54);
\draw (1.4,1.54)-- (1.13,1.57);
\draw (1.13,1.57)-- (1.18,2.07);
\draw (1.18,2.07)-- (0.95,2.46);
\draw (0.95,2.46)-- (0.65,2.66);
\draw (0.65,2.66)-- (0.88,2.93);
\draw (0.88,2.93)-- (0.44,3.07);
\draw (0.44,3.07)-- (0.35,3.46);
\draw (0.35,3.46)-- (0.29,3.73);
\draw (0.29,3.73)-- (-0.16,3.95);
\draw (-0.16,3.95)-- (0,4.57);
\draw (0,4.57)-- (-0.38,4.71);
\draw (-0.38,4.71)-- (-0.72,5.21);
\draw (-0.72,5.21)-- (-0.48,5.74);
\draw (-0.48,5.74)-- (-0.89,5.93);
\draw (-0.89,5.93)-- (-1.35,6.42);
\draw (-1.35,6.42)-- (-1,7);
\draw (-1,7)-- (-1.6,7.18);
\draw (-1.6,7.18)-- (-1.88,7.78);
\draw (-1.88,7.78)-- (-1.5,8);
\draw (-1.5,8)-- (-1.79,8.48);
\draw (-1.79,8.48)-- (-2.4,8.69);
\draw (-2.4,8.69)-- (-2.05,9.47);
\draw (-2.05,9.47)-- (-2.51,9.86);
\draw (0,7.34)--(1,7.34) node [right] {$B'$};
\draw (7.03,7)--(7.03,6) node [below] {$B$};
\draw (4.04,6.3) node [above left] {$\hat\g$};
\draw (5.77,7.15) node [below right] {$\g$};
\draw (-1.76,7.53) node [right] {$\g'$};
\fill (0.65,2.03) circle (1.5pt);
\draw (0.65,2.03) node [below] {$\bx$};
\end{tikzpicture}
    \caption{The paths $\g$ and $\g'$. The shifted and thickened version $\hat\g$ is the union of the shaded boxes. Due to its thickness it necessarily intersects $\g'$ and does so close to the intersection of the boxes $B$ and $B'$. The only $\bx$ yielding the intersection in the example is indicated by a dot.}
    \label{fig:DS2}
\end{figure}
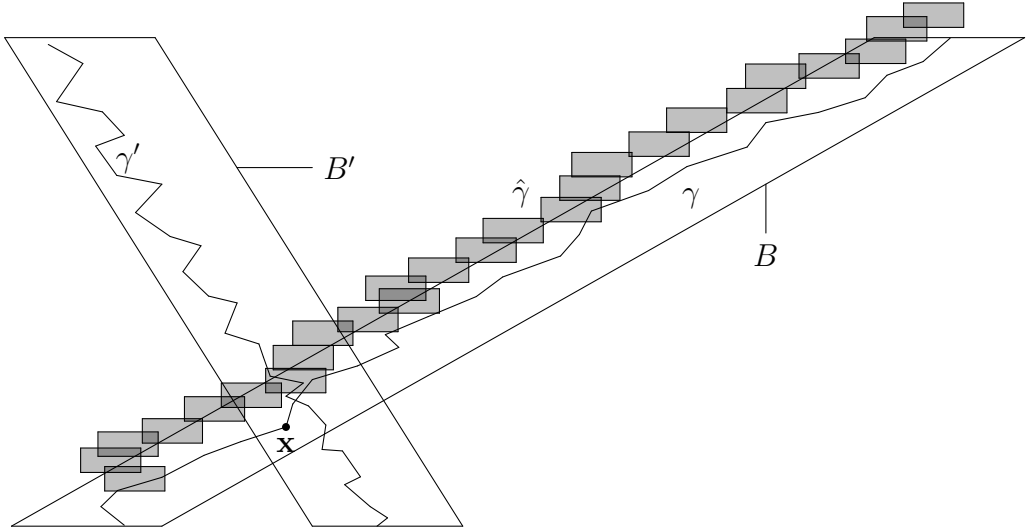

More precisely, using \cref{eq:crossing}, fix $\e>0$ and $\d>0$ small and $L$ large enough so that 
\[\bbP_p\left(S\to[B]S+(0,L)\right)>1-\d\]
with $B=B(\e L,L,\a)-2\e L \be_1$ and similarly for $B'=B(\e L,L,\b)+2\e L \be_1$ (see \cref{fig:DS2}). Fix two open paths $\g=(\ba_0,\ba_1,\dots,\ba_m)$ and $\g'=(\ba'_0,\ba'_1,\dots,\ba'_{m'})$ crossing $B$ and $B'$ respectively. Fix $v$ and large $n$ and $t$ as in \cref{eq:l} independent of all other constants. Let $\h$ be a set of additional infections, with each site at distance at most $O(t)$ from $B\cap B'$ infected independently with probability $\epsilon>0$. Then we claim that
\[\bbP_{p+\epsilon}\left(\ba_0\to[\g\cup\g'\cup\h]\ba'_{m'}\right)>\epsilon^{O(t^2)}.\]
To see this, simply consider the region
\[\hat\g=\bigcup_{\ba\in\g}(\ba+B_n+(\bv,t)),\]
which is a shifted and thickened version of $\g$ (see \cref{fig:DS2}). It is clear that $\hat\g\cap\g'\neq\varnothing$, so it suffices for one $\ba\in \g$ such that $(\ba+B_n+(\bv,t))\cap\g'\neq\varnothing$ to infect all possible sites for a time interval $t$ in $\h$. Since there are $O(t^2)$ of them, we obtain the desired result.

Hence, with positive probability we can go from $\ba_0$ to $\ba'_{m'}$ and, similarly, from $\ba'_0$ to $\ba_m$, which is just as good as having $\g\cap\g'\neq\varnothing$ for our purposes (except that the latter cannot be achieved by sprinkling). With this at hand, the approach of \cite{Durrett87} works without the annoying hypothesis (H3) to renormalise GOSP to 2dOP with parameter close to $1$. Consequently, $\a(p)>\b(p)$ does imply $\theta(p)>0$, but also, since the probability of a renormalised site being open is continuous in $p$ and arbitrarily close to $1$, \cref{th:DS} follows (see \cite{Durrett84} for more details).

\section*{Acknowledgements}
We would like to thank Barbara Dembin, Justin Salez, Cristina Toninelli and Daniel Valesin for helpful remarks.

\bibliographystyle{plain}
\bibliography{Bib}
\end{document}